\documentclass{amsart}

\usepackage{amssymb}
\usepackage{a4wide}

\usepackage{amsthm}
\usepackage{graphicx}
\usepackage[T2A]{fontenc}
\usepackage[latin1]{inputenc}

\renewcommand{\i}{{i'}}
\renewcommand{\j}{{j'}}
\renewcommand{\l}{{l'}}
\renewcommand{\k}{{k'}}

\newcommand{\R}{\mathbb R}

\newcommand{\bg}{\bar g}

\newcommand{\J}[4]{J_{\ \ #2}^{#1} J_{\ \ #4}^{#3}}
\newcommand{\Jij}{\mathcal J^{\i \j}_{i j}}
\newcommand{\Jkl}{\mathcal J^{\k \l}_{k l}}
\newcommand{\Jkj}{\mathcal J^{\k \j}_{k j}}
\newcommand{\const}{\mbox{\rm const}}
\newcommand{\Mu}{{\mathcal M}}
\newcommand{\Id}{{\bf 1}}

\newcommand{\weg}[1]{}
\newcommand{\spann}{\mathrm{span}}

\theoremstyle{plain}
\newtheorem{thm}{Theorem}
\newtheorem*{thm*}{Theorem}
\newtheorem{lem}{Lemma}

\newtheorem{cor}{Corollary}
\newtheorem*{con}{Convention}

\theoremstyle{definition}
\newtheorem{defn}{Definition}

\newtheorem{exmp}{Example}
\theoremstyle{remark}
\newtheorem{rem}{Remark}

\title{The only  Kähler manifold with degree of mobility $\geq 3$
  is $(\mathbb{C}P(\mbox{\rm n}), g_{\textrm{Fubini-Study}})$}
  \author{A. Fedorova, V. Kiosak, V.S. Matveev, S. Rosemann}
\address[A. Fedorova, V.S. Matveev, S. Rosemann]{Institute of Mathematics\\ FSU Jena\\ 07737 Jena, Germany}
\email{Aleksandra.Fedorova@uni-jena.de, vladimir.matveev@uni-jena.de, stefan.rosemann@uni-jena.de}
\address[V. Kiosak]{}
\email{vkiosak@ukr.net}

\thanks{partially supported by  GK 1523  and SPP 1154  of DFG}
\begin{document}

\begin{abstract}
  The degree of mobility of a (pseudo-Riemannian) Kähler metric is the
  dimension of the space of metrics $h$-projectively equivalent to
  it. We prove that a metric on a closed connected manifold can not have the
  degree of mobility $\ge 3$ unless it is essentially the Fubini-Study
  metric, or the $h$-projective equivalence is actually the affine
  equivalence. As the main application we prove an important special
  case of the classical conjecture attributed to Obata and Yano, 
  stating that a closed manifold admitting an essential group of
  {\it $h$-projective transformations}  is $(\mathbb{C}P(n),
  g_{Fubini-Study})$ (up to  multiplication of the
  metric by a constant). An additional result is the generalization of
  a certain result of Tanno 1978 for the pseudo-Riemannian situation.\\
  
\end{abstract}

\maketitle

\section{Introduction}
\label{sec:hpr}

\subsection{$h$-planar curves}

Let $(g, J)$ be a Kähler structure on a manifold $M^{2n}$. We allow
the metric $g$ to have arbitrary signature. A curve $\gamma:I\to
M^{2n}$ is called \emph{$h$-planar}, if there exist functions
$\alpha(t)$, $\beta(t)$ such that the following ODE holds:

\begin{gather}
  \label{eq:eq1}
  \nabla_{\dot{\gamma}}\dot{\gamma}=\alpha\dot{\gamma}+\beta
  J(\dot{\gamma}).
\end{gather}

Actually,  equation~(\ref{eq:eq1}) can be written as an ODE $
(\nabla_{\dot{\gamma}}\dot{\gamma} )\wedge \dot{\gamma} \wedge
J\dot{\gamma}=0$ on $\gamma$ only; but since this ODE is not in the
Euler form, there exist a lot of different $h$-planar curves with
the same initial data $\gamma(t_0), \dot\gamma(t_0)$. Nevertheless,
for every chosen functions $\alpha$ and $\beta$,  equation
(\ref{eq:eq1}) is an ODE of second order in the Euler form, and has
an unique solution with arbitrary initial values $\gamma(t_0),
\dot\gamma(t_0)$.

\par
Let us recall basic properties and basic examples of $h$-planar
curves.

\begin{exmp}  The property of a curve to be $h$-planar
survives after the reparametrization of the curve. In particular every 
(reparametrized) geodesic of $g$ is an $h$-planar curve. This is the reason why $h-$planar curves are also called {\it almost geodesics} or {\it complex geodesics} in the literature. 
\end{exmp}

\begin{exmp}
  \label{ex:1}
  Consider a $2$-dimensional Riemannian Kähler manifold, i.e. a Riemannian
  surface $(M^2, g)$ with the induced complex structure $J$. For this
  Kähler manifold every curve on $M^2$ is $h$-planar, since
  $\mbox{span}\{\dot \gamma(t) ,J(\dot \gamma(t))\}$ coincides with
  the whole $T_{\gamma(t)}M$ for $\dot\gamma(t)\ne 0$.
\end{exmp}

\begin{exmp}
  \label{ex:ex2}
  Consider $\mathbb{R}^{2n}=\mathbb{C}^{n}$ with the standard metric
  $g= \sum_{j=1}^n dz^j d\bar z^j$ and with the standard complex
  structure $J$ (acting by multiplication by  the imaginary unit
  $i$).\\
  Then, a curve $\gamma$ is $h$-planar if and only if it lies on a
  certain ``complex line'' $\mbox{Span}\{v,J(v)\}$ (for a certain
  $v\ne \vec 0$).
\end{exmp}

\begin{exmp}
  \label{ex:ex3}
  Consider the complex projective space
  \begin{align*}
    \mathbb{C}P (n)=\{1\mbox{-dimensional complex subspaces of }\mathbb{C}^{n+1}\}
  \end{align*}
with the standard complex structure $J= J_{standard}$.  The unitary
  group $U(n+1)$ acts naturally transitively by holomorphic
  transformations on $\mathbb{C}P(n)$. Since the group $U(n+1)$ is
  compact, there exists a Kähler metric on $\mathbb{C}P(n)$ invariant
  with respect to $U(n+1)$. This metric is unique up to multiplication
  by a constant and is called \emph{the Fubini-Study metric}, we
  denote it by the symbol $g_{FS}$. By an appropriate choice of the
  constant, $g_{FS}$ becomes a Riemannian metric of constant
  holomorphic sectional curvature equal to $1$ and we determine
  $g_{FS}$ uniquely by this choice.  Let $\pi$ be the standard
  projection $\pi: \mathbb{C}^{n+1}\setminus\{0\}\to
  \mathbb{C}P(n)$. We call a subset $L\subseteq\mathbb{C}P(n)$ a
  \emph{projective line}, if $L$ is the image of a $2$-dimensional
  complex subspace of $\mathbb{C}^{n+1}$ under the projection $\pi$.

  \par
  Let us see that every curve $\gamma$ lying on a certain projective
  line $L$ is $h$-planar (and vice versa). Indeed, $L$ is a totally
  geodesic 2-dimensional submanifold (for example since there exists
  an element $f\in U(n+1)$ such that $L$ is the set of fixed points of
  $f$). Since $L$ is $J-$invariant, $(L, g_{FS|L}, J_{|L})$ is a
  two-dimensional Kähler manifold (as in Example~\ref{ex:1}); in
  particular every curve on $(L, g_{FS|L}, J_{|L})$ is $h$-planar. Since
  the restriction of the connection of $g_{FS}$ to $L$ coincides with the
  connection of $g_{FS|L}$, every curve $h$-planar with respect to
  $(g_{FS|L}, J_{|L})$ is also $h$-planar with respect to $(g_{FS}, J)$. Now,
  every initial data $\gamma(0), \dot\gamma(0)$ and every functions
  $\alpha(t), \beta(t)$ can be realized by a $h$-planar curve lying on
  an appropriate projective line. Thus, a curve is $h$-planar if and
  only if it lies on a certain projective line $L$.
\end{exmp}

\subsection{ $h$-projectively equivalent metrics }

\begin{defn}[$h$-projectivity]
  Two metrics $g$ and $\bg$ that are Kähler with respect to the same
  complex structure $J$ are called \emph{$h$-projectively equivalent},
  if each $h$-planar curve of $g$ is an $h$-planar curve of $\bg$ and
  vice versa.
\end{defn}

\begin{exmp}
  If the metrics $g$ and $\bar g$ are Kähler with respect to the same
  complex structure $J$ and  are {\it affinely equivalent} (i.e., if their
  Levi-Civita connections $\Gamma$ and $\bar \Gamma$ coincide), then
  they are $h$-projectively equivalent. Indeed,  equation
 ~(\ref{eq:eq1}) for the first and for the second metric coincides if
  $\Gamma=\bar \Gamma$.
\end{exmp}

As we will see further, affine equivalence will be considered as a
special \emph{trivial} case of $h$-projectivity.

\begin{exmp}
  In particular, for every nondegenerate hermitian matrix $A=(a_{ij})
  \in Mat(n,n, \mathbb{C})$ the metric $\bar g = \sum_{i,j=1}^n
  a_{ij}dz^id\bar z^j $ is $h$-projectively equivalent to the metric 
  $g = \sum_{i=1}^n
  dz^id\bar z^i$
  from Example~\ref{ex:ex2}: indeed, the metric $\bar g$ is affinely
  equivalent to $g$ and is Kähler with respect to the same $J$. Though
  there exist other examples of metrics $h$-projectively equivalent to
  the metric from Example~\ref{ex:ex2}; they can be constructed
  similar to Example~\ref{ex:7}.
\end{exmp}

Let us now construct Kähler metrics $h$-projectively equivalent to the
Fubini-Study metric $g_{FS}$ on $\mathbb{C}P(n)$. The construction is
a generalization of the Beltrami's example of projectively equivalent
metrics, see~\cite{Beltrami}.

\begin{exmp}
  \label{ex:7}
  Consider a complex linear transformation of $\mathbb{C}^{n+1}$ given
  by a matrix $A\in GL_{n+1}(\mathbb{C})$ and the induced mapping
  $f_A:\mathbb{C}P^n \to \mathbb{C}P^n$ defined by $f_A(\pi(x))= \pi(Ax)$. Since
  the mapping $f_A$ preserves the complex lines $L$ and since by
  Example~\ref{ex:ex3} $h$-planar curves are those lying on a certain
  projective line $L$, the pullback $g_A:= f_A^*g_{FS}$ is
  $h$-projectively equivalent to $g_{FS}$. For further use let us note
  that the metric $g_A$ is isometric or affinely equivalent to
  $g_{FS}$ if and only if $A$ is proportional to a unitary matrix.
\end{exmp}

\subsection{PDE-system for $h$-projectively equivalent metrics and the
  degree of mobility. } 
\label{PDE-system}

Let $J$ be a complex structure on $M^{2n}$ and let $g$ and $\bg$ be
two metrics on $M^{2n}$ such that $(g, J)$ and $(\bg ,J)$ are Kähler
structures. We consider the following $(0,2)$-tensor $a_{ij}$ on $M$:

\begin{equation}
  \label{eq:a}
  a_{i j} = \left(\frac{\det{\bg}}{\det g}\right)^\frac{1}{2(n+1)}
  g_{i \alpha} \bg^{\alpha \beta} g_{\beta j},
\end{equation}
where $\bg^{\alpha \beta}$ is the $(2,0)$-tensor dual to $g_{\alpha
  \beta}$:
$\bg^{\alpha \beta}\bg_{\beta\gamma}=\delta^{\alpha}_{\gamma}.$\\
Obviously $a_{ij}$ is a hermitian, symmetric and non-degenerate
$(0,2)$-tensor.

\begin{con}
  We work in tensor notations. In particular we denote
  by ``comma'' the covariant differentiation with respect to the
  Levi-Civita connection defined by $g$, i.e., for example $T_{ij
    ,k}=\nabla_{k}T_{ij} $ for a (0,2)-tensor $T$. We sum with respect
  to repeating indices and use the metric $g$ to raise and lower indices,
  for example $J_{jk}=g_{j\alpha}J^{\alpha}_{ \ \ k}$ is the Kähler
  $2$-form corresponding to $g$. All indices range from $1$ to $2n$;
  the greek indices $\alpha,\beta,...$ also range from $1$ to $2n$ and
  will be mostly used as summation indexes (``dummy'' indices in
  jargon). We also introduce the following notation: for every
  $1-$form $\omega_i$ we denote by $\bar{\omega}_{i}=J_{\ \ i}^{
    \alpha}\omega_{\alpha}$ the ``multiplication'' of $\omega$ with
  the complex structure $J$.
\end{con}

The following statement  due to Mikes and Domashev plays an important role in the theory of
$h$-projectivity; it reformulates the condition ``$\bar g$ is
$h$-projectively equivalent to $g$'' to the PDE-language.\\

\begin{thm}[\cite{DomMik1977,DomMik1978}]
  \label{thm:mikes}
Let $(g,J)$ and $(\bar{g},J)$ be two Kähler structures on $M^{2n}$.
Then, $\bg$ is $h$-projectively equivalent to $g$ if and only if there
  exists a $(0,1)-$tensor $\lambda_i$ such that $a_{ij}$ given by
  \eqref{eq:a} satisfies
  \begin{equation}
    \label{eq:hpr}
    a_{i j,k}=\lambda_i g_{j k} + \lambda_j g_{i k} -
    \bar{\lambda}_{i}J_{jk}-\bar{\lambda}_{j}J_{ik}
  \end{equation}
\end{thm}

One can and should regard  equation~(\ref{eq:hpr}) as a PDE-system
on the unknown $(a_{ij}, \lambda_i)$ whose coefficients depend on the
metric $g$. Let us mention though that it is possible to consider
\eqref{eq:hpr} as a PDE-system on the unknown $(a_{ij})$ only: Indeed,
contracting \eqref{eq:hpr} with $g^{ij}$ we obtain
$\left(a_i^i\right)_{,k}= 4 \lambda_k$ (which in particular
implies that the covector $\lambda_i$ is a gradient, i.e.,
$\lambda_{i,j}=\lambda_{j,i}$).

\par
Note that the formula \eqref{eq:a} is invertible. Then, the set of the metrics $\bar g$ \
$h$-projectively equivalent to $g$ is essentially the same as the set
of the hermitian and symmetric solutions of \eqref{eq:hpr} (the only difference is the case
when $a_{ij}$ is degenerate; but since adding $\const\cdot g_{ij}$ to
$a_{ij}$ does not change the property of $a_{ij}$ to be a solution, this
difference is not important). Indeed, one can 
 show that if $(g,J)$ is Kähler, $a_{i j}$ is hermitian, symmetric, nondegenerate
and satisfies~\eqref{eq:hpr} for a certain $\lambda_i$, then the metric
$\bar g$ constructed via~\eqref{eq:a} is also Kähler  with respect to $J$.

\par
We see that the PDE-system~(\ref{eq:hpr}) is linear, hence the set of
its solutions is a linear vector space.

\begin{defn}
  The \emph{degree of mobility} of a Kähler metric $g$ is the
  dimension of the space of solutions $(a_{ij}, \lambda_i)$
  of \eqref{eq:hpr}, where $a_{ij}$ is symmetric and hermitian.
\end{defn}

\begin{rem}
  The degree of mobility $D$ is at least $1$ and is finite (assuming 
  dim$(M)\ge 4$; in the  two-dimensional case, 
  every two conformally equivalent metrics are $h$-projectively equivalent),
 $1\leq   D<\infty$. Indeed, $g$ itself is always a solution of~(\ref{eq:hpr})
  (with $\lambda_i\equiv 0$), implying $D\ge 1$. We will not make use
  of the fact that $D$ is finite, in fact $D\leq (n+1)^2$,  but it
  will be a direct consequence of Section~\ref{sec:proj} (and follows
  for example from~\cite[Theorem 2]{DomMik1978}).
\end{rem}

\begin{con}
  The equation~(\ref{eq:hpr}) plays a fundamental role in our
  paper. Whenever we speak about a  solution $(a_{ij}, \lambda_i)$
  of this equation, we assume that $a_{ij}$ is symmetric and
  hermitian. One of the reasons for it is that if $a_{ij}$ is
  constructed by \eqref{eq:a}, then it is automatically symmetric and
  hermitian. The second reason is that the procedure of symmetrization
  and hermitization
  $$T_{ij}\mapsto \frac{1}{4}T_{\alpha \beta} \left(
    \delta^\alpha_i \delta^\beta_j + \delta^\alpha_j \delta^\beta_{i}
    + J^\alpha_{\ \ i} J^\beta_{\ \ j}
    + J^\alpha_{\ \ j} J^\beta_{\ \ i} \right)$$
  does not affect the right-hand side of the equation; so without loss
  of generality we can always think that $a_{ij}$ in~(\ref{eq:hpr}) is
  symmetric and hermitian.
\end{con}

\begin{rem}
  \label{aff}
  For further use, let us note that if $\lambda_i\equiv 0$, then the
  metric $\bar g$ corresponding to $a_{ij}$ is affinely equivalent to
  $g$ (if it exists, i.e., if $a_{ij}$ is nondegenerate).
\end{rem}

\subsection{Main result}

Our main result is the following
\begin{thm}
  \label{thm:main}
  Let $(M^{2n},g,J)$ be a closed  connected Kähler manifold of degree of
  mobility $D\geq 3$ and of real dimension $2n\geq 4$. Then
  \begin{itemize}
  \item there is a constant $c\in \mathbb{R}$, $c\ne 0$,  such that $(M^{2n},c\cdot g,
    J)$ is  $(\mathbb{C}P(n),g_{FS}, J_{standard})$
    where $g_{FS}$ denotes the Fubini-Study metric on $\mathbb{C}P(n)$
    with the standard complex structure\\
    \\or\\
  \item each Kähler metric $\bg$, h-projectively equivalent to $g$, is
     affine equivalent to $g$.
  \end{itemize}
\end{thm}

In other words, a closed Kähler manifold $(M^{2n}, g, J)$ which is not
(a quotient of) $(\mathbb{C}P(n), \const\cdot g_{FS}, J_{standard})$
can not have $D\ge 3$ unless every metric $h$-projectively equivalent
to $g$ is affinely equivalent to $g$.

\par We would like to point out that we do not assume in
Theorem~\ref{thm:main} that the metric $g$ is Riemannian: an essential
part of the proof is to show that it must be definite (i.e., that
$\const \cdot g$ is Riemannian for an appropriate constant).

\subsection{ All conditions in Theorem~\ref{thm:main} are necessary. }

The assumption  $D\geq 3$ is necessary. Indeed, a construction of
a Kähler metric $g$ on $\mathbb{C}P(n)$ of non-constant holomorphic
sectional curvature such that it admits a metric $\bar g$ that is
$h$-projectively equivalent to $g$, but not affinely equivalent to $g$
can be extracted from~\cite{Kiyo1997}. In a certain sense, Kiyohara
found a way how one can perturb a pair of $h$-projectively equivalent
metrics on a closed manifold such that they remain $h$-projectively
equivalent. The space of perturbations is big and depends on
functional parameters. Perturbing $h$-projectively equivalent metrics
from Example~\ref{ex:7}, we obtain (for generic parameters of the
perturbation) metrics on $\mathbb{C}P(n)$ of non-constant holomorphic
sectional curvature admitting non-trivial $h$-projectivity. More examples can be extracted from \cite{ApostolovII}, see  discussion at the end of Section \ref{history}.

The assumption that the manifold is closed is also necessary. The
simplest examples of local metrics different from $g_{FS}$ with big
degree of mobility are due to~\cite{Tashiro1956}, see also
\cite{Fujimura,Tashibana}: it was shown that (locally) a metric of
constant holomorphic curvature (even if the metric is not positive
definite and the sign of the curvature is negative) admits a huge
space of $h$-projectively equivalent metrics. One can also construct
examples of (local) metrics of non-constant holomorphic curvature with
degree of mobility $\ge 3$ using the results of~\cite[\S2.2]{Mikes}.

\par
The second possibility in Theorem~\ref{thm:main} (when $g$ and $\bar
g$ are affinely equivalent) is also necessary. Indeed, consider the
direct product of three Kähler manifolds
$$(M_1, g_1, J_1) \times (M_2, g_2, J_2) \times (M_3, g_3, J_3).$$
It is a Kähler manifold diffeomorphic to the product $M_1\times M_2
\times M_3$, the metric is the sum of the metrics $g_1 + g_2 +g_3$,
and the complex structure is the sum of the complex structures. Then,
for any constants $c_1, c_2, c_3 \ne 0$, the metrics $c_1 \cdot g_1 +
c_2 \cdot g_2 + c_2 \cdot g_3$ is $h$-projectively equivalent to $g_1
+ g_2 +g_3$ (because they are  affinely equivalent to it), i.e., the
degree of mobility of $g_1 + g_2 +g_3$ is at least 3. If $M_i$ are
closed, then $M_1\times M_2 \times M_3$ is closed as well. Of
course, the metric $g_1 + g_2 +g_3$ is not $\const\cdot g_{FS}$.

\subsection{History, motivation, and first applications}

\subsubsection{History and motivation.}
\label{history}

$h$-planar curves and $h$-projectivity of Kähler metrics where
introduced in~\cite[\S\S 9-10]{Otsuki1954}. Otsuki and Tashiro did not
explain explicitly their motivation, from the context one may suppose
that they tried to study projectively equivalent metrics (the
definition is in Section~\ref{relation}) in the Kähler situation,
found out that they are not interesting (impossible except of few
trivial examples), and suggested a Kähler analog of projectively
equivalent metrics. Actually, it was one of the main trend of their
time to adapt Riemannian objects to the Kähler situation, see for
example the book~\cite{Yanobook} (where many objects were generalized
to the Kähler situation; $h-$projectively equivalent metrics are in
the last chapter of this book).

\par The notion turned out to be interesting and successful, there are
a lot of papers studying $h$-projectivity and its generalizations, see
for example the recent survey~\cite{Mikes}.  At a certain period of
time $h-$projectivity was one of the main research topics of the
Japanese and Soviet (mostly Odessa and Kazan) geometry schools.  At
least two books,~\cite{Sinjukov} and~\cite{Yanobook}, have chapters on
$h-$projectively equivalent metrics.

\par One of the mainstreams in the theory of $h$-projectivity is to
understand the group of {\it $h$-projective transformations}, i.e., the
group of diffeomorphisms of $(M^{2n}, g, J)$ that preserve the complex
structure and send the metric to a metric that is $h$-projectively
equivalent to $g$. This set is obviously a group, Ishihara
\cite{Ishihara1957} and {Yoshimatsu}~\cite{Yoshimatsu} have shown that
it is a finite dimensional Lie group and the challenge was to
understand the codimension of the group of affine transformations or
isometries in this group, see for example~\cite{Ishihara1957,
  HiramatuK,Yano1981,Akbar,HF,Mikes}.

\par As it follows from Example~\ref{ex:7}, the group of
$h$-projective transformations of $(\mathbb{C}P(n),g_{FS},
J_{standard})$ is much bigger than its subgroup of affine
transformations. A classical conjecture (in folklore this conjecture
is attributed to Obata and Yano, though we did not find a reference
where they formulate it explicitly) says that, {\it on closed  Riemannian
Kähler manifolds that are not $(\mathbb{C}P(n), \const \cdot
g_{FS} ,J_{standard})$, the connected component of the group of
$h$-projective transformations contains isometries only.}  In
particular, in the above mentioned papers~\cite{Ishihara1957,
  Ishihara1961, HiramatuK,Yano1981,Akbar}, the conjecture was proved
under certain additional assumptions; for example, the additional
assumption in~\cite{HiramatuK,Yano1981,Akbar} was that the scalar
curvature of the metric is constant. 

For the Riemannian metrics, the Yano-Obata conjecture was proved in the recent paper \cite{yano-obata}. 
The proof uses different techniques to those employed in the present article, but does rely in part on certain results (Theorem \ref{thm:main} and Section \ref{subsec:globalB}) of the present paper.

 In Section~\ref{applications},  we give new results assuming that the metric has arbitrary signature. In particular, we  show that the codimension of the
subgroup of isometries in the group of $h$-projective transformation
is at most one.

\par Recent interest to $h$-projectivity is in particular due to an
unexpected connection between $h$-projectively equivalent metrics and
integrable geodesic flows: it appears that the existence of $\bar g$
$h$-projectively equivalent to $g$ allows to construct quadratic and
linear integrals for the geodesic flow of $g$, see for example
\cite{Top2003,Kiyohara2010}. Theorem \ref{thm:main} shows that there is no metric (except of Fubini-Study) on a closed Kähler manifold such that its geodesic flow is superintegrable with integrals coming from $h$-projectively equivalent metrics.  
 
 Additional interests to $h$-projective equivalence is due to its connection with the so called {\it hamiltonian $2$-forms} defined  and investigated in Apostolov et al  \cite{ApostolovI,ApostolovII,ApostolovIII,ApostolovIV}. It is easy to see that  a  hamiltonian $2$-form   is essentially the same as a $h$-projectively equivalent 
metric $\bar g$, since   the defining equation \cite[equation $(12)$]{ApostolovI} of a hamiltonian $2$-form is algebraically equivalent to the equation \eqref{eq:hpr} from Theorem \ref{thm:mikes}. The motivation of Apostolov et al  to study hamiltonian $2$-forms is different from that of Otsuki and Tashiro and is  explained in \cite{ApostolovI,ApostolovII}. Roughly speaking, they observed that many interesting problems on  Kähler manifolds lead to  hamiltonian $2$-forms and suggested to study them. The motivation is justified in \cite{ApostolovIII,ApostolovIV}, where they indeed constructed  new interesting and useful examples of Kähler manifolds. There is also a direct connection between $h$-projectively equivalent metrics and  conformal Killing  (or twistor) $2$-forms  studied in \cite{Moroianu,Semmelmann0,Semmelmann1}, see Appendix A of \cite{ApostolovI} for details. 

In private communications with the authors of \cite{ApostolovI,ApostolovII,ApostolovIII,ApostolovIV} we got informed that they did not know that the object they considered  was studied before under another name. Indeed, they re-derived certain facts that were well known in the theory of $h$-projectively equivalent metrics. 
On the other hand, the papers    \cite{ApostolovI,ApostolovII,ApostolovIII,ApostolovIV} contain several solutions of the problems studied in the framework of $h$-projectively equivalent metrics; in particular they gave a global description of metrics admitting hamiltonian 2-forms providing us with new nontrivial examples of $h$-projectively equivalent metrics.

 \subsubsection{First applications: special case of the Yano-Obata
   conjecture.}
 \label{applications}

 Let $(M^{2n}, g, J)$ be a Kähler manifold. Recall that a
 diffeomorphism $f:M\rightarrow M$ is called a \emph{$h$-projective
   transformation} if it preserves the complex structure $J$ and sends
 the metric $g$ to a metric that is $h$-projectively equivalent to
 $g$. The set of all $h$-projective transformations of $(M^{2n}, g, J)$
 forms a Lie group which we denote by $\mbox{HProj}$. We denote by
 $\mbox{HProj}_{0}$ its connected component containing the
 identity. The groups of affine transformations and isometries of $M$
 preserving the complex structure and their connected components
 containing the identity will be denoted by $\mbox{Aff}(g,J)$,
 $\mbox{Iso}(g,J)$, $\mbox{Aff}_{0}(g,J)$, and $\mbox{Iso}_0(g,J)$,
 respectively.

 \begin{cor}
   \label{cor1}
   Let $(M^{2n},g,J)$ be a closed connected Kähler
   manifold of dimension $2n \geq 4$. Assume that for every $\const\ne 0$ the manifold
   $(M^{2n},g,J)$ is not  $(\mathbb{C}P(n),\const\cdot
   g_{FS}, J_{standard})$. Then the group $\mbox{Iso}_{0}(g,J)$ has
   the codimension at most one in the group $\mbox{HProj}_{0}$, or
   $\mbox{HProj}= \mbox{Aff}(g,J)$.
 \end{cor}

\begin{proof}
  First assume $D=1$. This means that each metric that is
  $h$-projectively equivalent to $g$, is proportional to it. Thus,
  every $h$-projective transformation is a homothety. Since the
  manifold is closed, every homothety is an isometry implying $\mbox{HProj}=
  \mbox{Iso}(g,J)$.

  \par
  Assume now $D\geq 3$. Then, by Theorem~\ref{thm:main}, every $\bar
  g$ $h$-projectively equivalent to $g$ is affinely equivalent to $g$
  implying $\mbox{HProj}= \mbox{Aff}(g,J)$.

  \par
  The remaining case is $D=2$. We need to show that the Lie-algebra of
  $\mbox{Iso}_{0}(g,J)$ has codimension at most one in the Lie-algebra
  of $\mbox{HProj}_{0}$.

  \par
  Let $u,v$ be infinitesimal $h$-projective transformations,
  i.e. vector fields on $M$ generating $1$-parameter groups of 
 $h$-projective transformations. We need to show that their
  certain linear combination is a Killing vector field. Let us first
  construct a mapping $\Psi: u\mapsto a_u$ sending an infinitesimal
  $h$-projective transformation to a solution of \eqref{eq:hpr}.

  \par
  We denote by $\Phi^{u}_{t}$ the  flow of $u$ and define
  $g_{t}:=(\Phi^{u}_{t})^{*}g$. As we recalled in  Section \ref{PDE-system} (see Theorem \ref{thm:mikes} there), the
  $(0,2)-$tensor $a(t)_{ij}$ given by (in matrix notation)
  \begin{gather}
    a(t)=\left(\frac{\mbox{det}g_{t}}{\mbox{det}g}\right)^{\frac{1}{2(n+1)}}gg_{t}^{-1}g\nonumber
  \end{gather}
  satisfies equation~(\ref{eq:hpr}). Taking the derivative at $t=0$,
  and replacing the $t-$derivatives of tensors by Lie derivatives, we
  obtain that the $(0,2)-$tensor
  \begin{gather}
    a_{u}:=L_{u}g-\frac{\mbox{trace}\,g^{-1}L_{u}g}{2(n+1)}g\nonumber
  \end{gather}
  satisfies equation~(\ref{eq:hpr}).

  \par
  We define then the mapping $\Psi$ by $\Psi(u)= a_u$. The mapping is
  clearly linear in $u$. Since the two-dimensional space of the
  solutions of~(\ref{eq:hpr}) contains the one-dimensional subspace
  $\{c\cdot g\mid c\in \mathbb{R}\}$, for every two infinitesimal
  $h$-projective transformations $u, v$ there exists a linear
  combination $b u + d v $ such that $\Psi(b u + dv) = c g$ (for a
  certain $c\in \mathbb{R}$). Let us show that $b u + d v $ is a
  Killing vector field. We have:
  \begin{equation}
    \label{-1} L_{bu +dv}g-\frac{\mbox{trace}\,g^{-1}L_{bu+
        dv}g}{2(n+1)}g = cg.
  \end{equation}
  Multiplying this (matrix) equation by the inverse matrix of $g$ and
  taking the trace, we obtain
  $$\mbox{trace}(g^{-1}L_{bu+ dv}g) - \tfrac{2n}{2(n+1)} \mbox{trace}(g^{-1}L_{bu+ dv}g) = 2n c.$$
  Thus, $\mbox{trace}(g^{-1}L_{bu+ dv}g) = 2n(n+1)c$. Substituting
  this in \eqref{-1}, we obtain that $L_{bu+ dv}g = c(1-n)\cdot
  g$. Then, $bu+ dv$ is an infinitesimal homothety. Since the manifold
  is closed, any infinitesimal homothety is a Killing vector field
  implying that $bu+ dv $ is a Killing vector field as we claimed.
\end{proof}

\subsection{Additional motivation: new methods for the investigation of the global behavior of $h$-projectively equivalent pseudo-Riemannian metrics.}

In many cases, local statements about Riemannian metrics could be
generalised for the pseudo-Riemannian setting, though sometimes this
generalisation is difficult. As a rule, it is very difficult to
generalize global statements about Riemannian metrics to the
pseudo-Riemannian setting. The theory of $h$-projectively equivalent
metrics is not an exception: certain local results could be
generalized without essential difficulties.  Up to now, no global (say
if the manifold is closed) methods for the investigation of
$h$-projectively equivalent metrics were generalized for the
pseudo-Riemannian setting.

\par More precisely, virtually every global result (see for example
the surveys~\cite{Mikes,Sin2}) on $h-$projectively equivalent Riemannian
metrics was obtained by using the so-called ``Bochner technique'',
which requires that the metric is positively defined.

\par Our proofs (we explain the scheme in Section \ref{plan}) use
essentially new methods (in Section~\ref{relation} we explain that
these methods were motivated by new results in the theory of
projectively equivalent metrics).  We expect further applications of
these new methods in the theory of $h$-projectively equivalent
metrics, and in other parts of differential geometry.

\subsection{Additional result: Tanno-Theorem for pseudo-Riemannian
  metrics} \label{1.8} 
Let us recall the following classical result of Tanno and Hiramatu:

\begin{thm}[\cite{Tanno1978},\cite{HiramatuK}]
{ Let  $f$ be  a non-constant smooth function on a closed   Riemannian Kähler manifold
  $(M^{2n},g,J)$ of dimension $2n\ge 4$ such that the equation
  \begin{gather}
    \label{eq:tanno}
    f_{,ijk}=\kappa(2f_{,k} \cdot g_{ij}+ f_{,i} \cdot  g_{j k} +
    f_{,j} \cdot g_{i k}-\bar{f}_{,i} \cdot J_{jk}-\bar{f}_{,j}\cdot J_{ik}).
  \end{gather}
  is  fulfilled (for a certain constant $\kappa$).  Then,  $\kappa< 0$  and  $(M^{2n},g,J)$ has
  constant holomorphic sectional curvature $-{4\kappa}$.  In particular, $(M^{2n}, -4 \kappa \cdot g, J)$ is  $(\mathbb{C}P(n), g_{FS}, J_{standard})$. } \end{thm}

  More precisely, Tanno \cite{Tanno1978} assumed that $\kappa<0$;  in this case it is sufficient to require that the manifold is complete. Hiramatu \cite{HiramatuK} proved  that the equation can not have nonconstant solutions for $\kappa\ge  0$, if   the manifold is closed. One can construct  counterexamples to the  latter statement, if the manifold  is merely complete.

We will show  in Section \ref{newtan} that a part of the proof of our main result  gives also a
proof of the pseudo-Riemannian version of the above statement:
\begin{thm}
  \label{thm:tanno}
  Let  $f$ be  a non-constant smooth function on a closed  connected pseudo-Riemannian Kähler manifold
  $(M^{2n},g,J)$ of dimension $2n\ge 4$ such that the equation
 \eqref{eq:tanno} 
  is fulfilled (for a certain constant $\kappa$). Then, $\kappa\ne 0$ and 
   $(M^{2n}, -4 \kappa \cdot g, J)$ is  $(\mathbb{C}P(n), g_{FS}, J_{standard})$.
\end{thm}

\subsection{Plan of the proof.} \label{plan} 

We assume that $(M^{2n}, g, J)$ is a closed connected Kähler manifold
of dimension $2n\ge 4$. We divide the proof of Theorem~\ref{thm:main}
in four steps.

\begin{itemize}

\item In Section~\ref{sec:extsys}, assuming $D\geq 3$, we show that
  for every solution $(a_{ij},\lambda_i)$ of equation~(\ref{eq:hpr}) 
there exists a  constant $B\in\mathbb{R}$ and a function $\mu$
  such that the following ``extended'' system
  \begin{gather}
    a_{i j,k}=\lambda_i  g_{j k} + \lambda_j g_{i k} - \bar{\lambda}_{i}J_{jk}-\bar{\lambda}_{j}J_{ik}\nonumber\\
    \lambda_{i,j}=\mu g_{i j}+Ba_{i j}\nonumber\\
    \mu_{,i}=2B\lambda_i\nonumber
  \end{gather}
  is satisfied (see Theorem~\ref{thm:ext}).  For   $a_{ij}\neq\mbox{const}\cdot g_{ij}$, the constant
   $B$ is  uniquely determined by the metric (Corollary~\ref{cor:Bind}), i.e.,
  is the same for all solutions of~\eqref{eq:hpr} that are not proportional to $g$.

In Sections ~\ref{sec:B=0}, \ref{sec:proj}, \ref{sec:tan}  we will work with the above ``extended'' system
  only, i.e., we will  not use that the degree of mobility of $g$ is $\ge
  3$ anymore. We show that \emph{the existence of a solution $(a_{ij},
    \lambda_i, \mu)$ with $\lambda_i\not\equiv 0$ on a closed connected Kähler manifold
    implies that the metric is proportional to the Fubini-Study
    metric. } We proceed as follows:
       \item In Section~\ref{sec:B=0} (see Theorem~\ref{thm:B=0}),
  we show that $B\ne 0$ unless $\lambda_i\equiv 0$.
\item If $B\ne 0$, by replacing $g$ with $-B\cdot g$, without loss of generality we can assume
    $B=-1$. In Section~\ref{sec:proj},
  we show that, for $B= -1$, the metric $g$ is positively definite.
\item In Section~\ref{sec:tan}, we combine the results of the previous
  sections and  the  result of Tanno~\cite{Tanno1978} we recalled in Section \ref{1.8}, to show
  that our manifold is  $(\mathbb{C}P(n), \const\cdot
  g_{FS}, J_{standard})$. This concludes the proof of Theorem~\ref{thm:main}.
\end{itemize}

\subsection{ Relation with projective equivalence.}
\label{relation}

Two metrics $g$ and $\bar g$ on the same manifold are
\emph{projectively equivalent}, if every geodesic of $g$, after an
appropriate reparametrization, is a geodesic of $\bar g$. As we
already mentioned in Section~\ref{history}, we think that the notion
``$h$-projective equivalence'' was introduced as an attempt to adapt
the notion ``projective equivalence'' to Kähler metrics. It is
therefore not a surprise that certain methods from the theory of
projectively equivalent metrics could be adapted for  the $h$-projective
questions. For example, the above mentioned papers
\cite{HiramatuK,Yano1981,Akbar} are actually an $h$-projective analog
of the papers~\cite{Yamauchi1974,Hiramatu} (dealing with projective
transformations), see also~\cite{Hasegawa,Takeda}.  Moreover,
\cite{Yoshimatsu,Tashiro1956} are ``Kählerizations'' of
\cite{Ishihara1961,Tanaka}, and many results listed in 
the survey \cite{Mikes} are
``Kählerizations'' of those listed in \cite{Mikes1996}.

 The Yano-Obata conjecture is also an $h$-projective analog of the
so-called projective Lichnerowicz-Obata conjecture (recently proved in
\cite{Matveev2007,CMH}, see also~\cite{M2004,M2004bis}). There also
exists a conformal analog of this conjecture (the so called conformal
Lichnerowicz-Obata conjecture proved in \cite{Alekseevskii1972, Obata, Schoen}),
whose $CR-$analog was proved in~\cite{Schoen}, and finsler analog in \cite{MRTZ}.

\par We also used certain ideas from the theory of projectively
equivalent metrics. In particular, the scheme of the first part of the
proof of Theorem \ref{thm:main} is close to the scheme of the proof of
\cite[Theorem 1]{KioMat2010}, see also~\cite{M2006}, the scheme of the
second part of the proof is close to the proof of~\cite[Theorem
1]{Mounoud2010} (though the proofs in the present paper are
technically much more complicated than the proofs in
\cite{KioMat2010,Mounoud2010}).

\par Let us also recall that recently new methods for the
investigation of projectively equivalent metrics were suggested. A
group of these new methods came from the theory of integrable systems
and from the dynamical systems~\cite{KM2006,MT98,M2003,Inventiones2003}.  We
expect that these methods could also be adapted for  the investigation of
$h$-projectively equivalent metrics (first steps were already done in
\cite{Kiyohara2010}).
 Another group of new methods came from the geometric theory of
ODEs, see for example~\cite{BMM2008,M2010,BDE}. We expect that these
methods could also be adapted for  $h$-projective transformations.

Let  us also recall that  equation \eqref{eq:tanno} was introduced  in \cite{Tanno1978} as ``Kählerization''
 of   $  f_{,ijk}=\kappa(2f_{,k} \cdot g_{ij}+ f_{,i} \cdot g_{j k} +
    f_{,j} \cdot g_{i k})$.  The latter
    equation  appeared  independently  and was helpful 
    in many  parts of differential geometry:  in spectral geometry \cite{Tanno1978,Gallot1979}, in cone geometry \cite{Gallot1979,Cortes2009},         and in conformal and projective geometry (see \cite{Hiramatu,Tanno1978} and \cite{Matveev2010,Mounoud2010} for references). We expect that  equation  \eqref{eq:tanno} will be helpful in the  ``Kählerizations''  of these geometries.

\section{Local theory and extended system}
\label{sec:extsys}
The goal of Section~\ref{sec:extsys} is to prove the following
\begin{thm}
  \label{thm:ext}
  Let $(M^{2n},J,g)$ be a connected Kähler manifold of dimension $2n\geq
  4$. If the degree of mobility $D$ of $g$ is $\geq 3$, then for
  every solution $(a_{ij},\lambda_i)$ of~(\ref{eq:hpr}), such that $a_{ij}\neq \mbox{const}\cdot g$, there exists
  a unique constant $B$ and a scalar function $\mu$, such that the
  extended system
  \begin{equation}
    \label{eq:mg+Ba}
    \begin{array}{ll}
      a_{i j,k}&=\lambda_i g_{j k} + \lambda_j g_{i k} - \bar{\lambda}_{i}J_{jk}-\bar{\lambda}_{j}J_{ik}\\
      \lambda_{i,j}&=\mu g_{i j}+Ba_{i j}\\
      \mu_{,i}&=2B\lambda_i
    \end{array} \end{equation}
  is satisfied.
\end{thm}
We see that the first equation of~(\ref{eq:mg+Ba}) is precisely the
equation~(\ref{eq:hpr}), i.e., is fulfilled by assumptions. We would
like to note here that the second and the third equations are
\emph{not} differential consequences of the first one: they require
the assumption that the degree of mobility is $\ge 3$.

\par
The proof of the second equation is the lengthiest and trickiest part
of the proof of Theorem~\ref{thm:ext}.  After recalling  basic properties of $\lambda_{i}$ in Section \ref{sec:int}, we will first prove a pure
algebraic result (Lemma~\ref{lem:ext}). Together with Lemma \ref{lem:globind}, it will imply that the
equation $\lambda_{i,j}=\mu g_{i j}+Ba_{i j}$ holds in a neighborhood of
almost every point of $M$ for a certain \emph{function} $B$. Then, in
Lemma~\ref{lem:ext2} we show that, locally, in a neighborhood of
almost every point, the function $B$ is actually a constant.  The
constant $B$ and the function $\mu$ could a~priori depend on a
neighborhood of the manifold, the last step will be to show that $B$
and $\mu$ are the same for each neighborhood and, hence, are globally
defined (Section~\ref{subsec:globalB}). Now, the third equation of
Theorem~\ref{thm:ext} will be obtained as a differential corollary of
the first two.

Note also that $(g_{ij}, 0)$ is also a solution of    \eqref{eq:mg+Ba}, with $\mu=-B$, so Theorem \ref{thm:ext} 
holds for this solution except for the constant $B$ is not unique anymore. In Section \ref{sec:proj} we will consider 
$(g_{ij}, 0)$ as a solution of \eqref{eq:mg+Ba} with $B=-1$ and $\mu=1$. 

\subsection{Killing vector field for the geodesic flow of $g$}
\label{sec:int}
In this section we show that the $1$-form $\bar{\lambda}_{i}$ satisfies the Killing equation, a fact which we shall use several times during our paper. 
\begin{lem}[Folklore]
 \label{lem:skew}
Let $(M^{2n}, g, J)$ be a Kähler manifold of dimension $2n\ge 4$ and let $(a_{ij},\lambda_{i})$ be a solution of equation (\ref{eq:hpr}). Then
 $J$ anticommutes with $g_{ij}$, $a_{ij}$ and $\lambda_{i,j}$:
 \begin{gather}
   \label{eq:gJ} J_{\ \ i}^\alpha g_{\alpha j}=-g_{i \alpha} J_{\ \ j}^\alpha,\nonumber \\
   J_{\ \ i}^\alpha a_{\alpha j}=-a_{i \alpha} J_{\ \ j}^\alpha,\nonumber \\
   J_{\ \ i}^\alpha \lambda_{\alpha, j}=-\lambda_{i, \alpha} J_{\ \ j}^\alpha.\nonumber
 \end{gather}
\end{lem}
\begin{proof}
  The first equality is a part of the definition of Kähler metrics,
  the second property follows  from  our convention from
  Section~\ref{PDE-system}. The third equality is also somehow known: it follows immediately from  \cite[equation (13)]{DomMik1978} and \cite[Proposition 3]{ApostolovI}. For the convenience of the reader, we give its proof but it does not pretend to be new. 

  \par
  Differentiating~(\ref{eq:hpr}), we obtain
  \begin{align}
    a_{ij,kl}=\lambda_{i,l}g_{jk}+\lambda_{j,l}g_{ik}-\bar{\lambda}_{i,l}J_{jk}-\bar{\lambda}_{j,l}J_{ik}.\nonumber
  \end{align}
  Substituting this into the formula
  $a_{ij,kl}-a_{ij,lk}=R^{r}_{ikl}a_{rj}+R^{r}_{jkl}a_{ir}$ (which is
  fulfilled for every $(0,2)$-tensor $a_{ij}$) we obtain
  \begin{equation}
    \label{integr}
    \begin{array}{ll}
      a_{ij,kl}-a_{ij,lk}&=\lambda_{i,l}g_{jk}-\lambda_{i,k}g_{jl}+\lambda_{j,l}g_{ik}-\lambda_{j,k}g_{il}
      -\bar{\lambda}_{i,l}J_{jk}+\bar{\lambda}_{i,k}J_{jl}-\bar{\lambda}_{j,l}J_{ik}+\bar{\lambda}_{j,k}J_{il}\\
      &=R^{r}_{ikl}a_{rj}+R^{r}_{jkl}a_{ir}.
    \end{array}
  \end{equation}
  Multiplying this equation with $g^{jk}$ and summing with respect to
  repeating indices, we obtain:
  \begin{align}
    \label{eq:in}
    2n\lambda_{i,l}-\lambda_{i,l}+\lambda_{i,l}-g^{jk}\lambda_{j,k}g_{il}-
    0+\bar{\lambda}_{i,k}J^{k}_{\ \ l}-\lambda_{i,l}+g^{jk}\bar{\lambda}_{j,k}J_{il}\nonumber\\
    =(2n-1)\lambda_{i,l}-g^{jk}\lambda_{j,k}g_{il}+\bar{\lambda}_{i,k}J^{k}_{\ \ l}
    +g^{jk}\bar{\lambda}_{j,k}J_{il}=g^{jk}R^{r}_{ikl}a_{rj}+g^{jk}R^{r}_{jkl}a_{ir}.
  \end{align}
  Recall that $g_{ij}$ and $a_{ij}$ are hermitian and the curvature
  satisfies the symmetry relations
  \begin{align}
    R^{i}_{\alpha kl}J^{\alpha}_{\ \ j}=
    J^{i}_{\ \ \alpha}R^{\alpha}_{jkl}\mbox{ and }R^{i}_{j\alpha\beta}
    J^{ \alpha}_{\ \ k}J^{\beta}_{\ \ l}=R^{i}_{jkl}.\nonumber
  \end{align}
  Now, let us rename $i\rightarrow i'$ and $l\rightarrow l'$, multiply equation~(\ref{eq:in}) by $J^{i'}_{\ \ i} J^{l'}_{\ \ l}$, and sum
  with respect to repeating indices. We want to show that this operation
  does not change the right-hand side of the equation.  First
  we consider the second term on the right-hand side:
  \begin{align*}
    g^{jk}R^{r}_{jkl'}a_{i'r}J^{i'}_{\ \ i}J^{l'}_{\ \ l}=
    -g^{jk}R^{r}_{jkl'}a_{i'i}J^{i'}_{\ \ r}J^{l'}_{\ \ l}=
    -g^{jk}R^{i'}_{rkl'}a_{i'i}J^{r}_{\ \ j}J^{l'}_{\ \ l}=
    -g^{jk}R^{r}_{i'kl'}a_{ri}J^{i'}_{\ \ j}J^{l'}_{\ \ l}\\
    =g^{ji'}R^{r}_{i'kl'}a_{ri}J^{k}_{\ \ j}J^{l'}_{\ \ l}=
    g^{ji'}R^{r}_{i'jl}a_{ri}=g^{jk}R^{r}_{jkl}a_{ir}
  \end{align*}
  We see that this term remains unchanged. Similarly, for the first term
  on the right-hand side we have
  \begin{align*}
    g^{jk}R^{r}_{i'kl'}a_{rj}J^{i'}_{\ \ i}J^{l'}_{\ \ l}=
    g^{jk}R^{i'}_{ikl'}a_{rj}J^{r}_{\ \ i'}J^{l'}_{\ \ l}=
    -g^{jk}R^{i'}_{ikl'}a_{ri'}J^{r}_{\ \ j}J^{l'}_{\ \ l}\\
    =g^{jr}R^{i'}_{ikl'}a_{ri'}J^{k}_{\ \ j}J^{l'}_{\ \ l}=
    g^{jr}R^{i'}_{ijl}a_{ri'}=g^{jk}R^{r}_{ikl}a_{jr},
  \end{align*}
  which again shows that the operation above does not change this
  term. Thus, the right-hand side of~(\ref{eq:in}) remains unchanged,
  so the difference of the left-hand side of~(\ref{eq:in}) and the
  transformed left-hand side of~(\ref{eq:in}) must be zero. We obtain:
  \begin{align*}
    0=(2n-1)\lambda_{i',l'}J^{i'}_{\ \ i}J^{l'}_{\ \ l}-
    g^{jk}\lambda_{j,k}g_{i'l'}J^{i'}_{\ \ i}J^{l'}_{\ \ l}+
    \bar{\lambda}_{i',k}J^{k}_{\ \ l'}J^{i'}_{\ \ i}J^{l'}_{\ \ l}+
    g^{jk}\bar{\lambda}_{j,k}J_{i'l'}J^{i'}_{\ \ i}J^{l'}_{\ \ l}\\
    -(2n-1)\lambda_{i,l}+g^{jk}\lambda_{j,k}g_{il}-\bar{\lambda}_{i,k}
    J^{k}_{\ \ l}-g^{jk}\bar{\lambda}_{j,k}J_{il}\\
    =(2n-1)\bar{\lambda}_{i,k}J^{k}_{\ \ l}-(2n-1)\lambda_{i,l}-
    \bar{\lambda}_{i,k}J^{k}_{\ \ l}-\bar{\lambda}_{i',l}J^{i'}_{\ \ i}
    =(2n-2)(\bar{\lambda}_{i,k}J^{k}_{\ \ l}-\lambda_{i,l})
  \end{align*}
  Hence, $\bar{\lambda}_{i,k}J^{k}_{\ \ l}=\lambda_{i,l}$. Multiplying
  by $J^{l}_{\ \ j}$ and using that $\lambda_{i,j}$ is symmetric
  yields the desired formula
  $-\bar{\lambda}_{i,j}=\lambda_{i,l}J^{l}_{\ \
    j}=\lambda_{l,i}J^{l}_{\ \ j}=\bar{\lambda}_{j,i}$.
\end{proof}

\begin{cor}[\cite{ApostolovI}]
  \label{cor:killing}
  Let $(M^{2n},g,J)$ be a Kähler manifold of dimension $2n\geq 4$. If $(a_{ij}, \lambda_i) $ is a solution of equation~(\ref{eq:hpr}),
  then $\bar \lambda^i:= g^{i\alpha}\bar{\lambda}_{\alpha}$ is a
  Killing vector field for $g$.
\end{cor}

\begin{proof}
  A vector field $v^i$ is Killing, if and only if the Killing equation
  $v_{i,j} + v_{j,i}=0$ is satisfied. For the vector field $\bar
  \lambda^i$, the Killing equation reads $\bar \lambda_{i,j}+ \bar
  \lambda_{j,i}=0$ and is equivalent to the third equality of
  Lemma~\ref{lem:skew}.
\end{proof}

\begin{cor}
  \label{cor:killing2}
  Let $(a_{ij}, \lambda_i) $ be a solution of equation~(\ref{eq:hpr})
  on a connected Kähler manifold $(M^{2n},g,J)$ of dimension $2n\geq 4$. If $\lambda_i\ne 0$ at a point, then
  $\lambda_i\ne 0$ at almost every point.
\end{cor}

\begin{con}
  Within the whole paper we understand ``almost everywhere'' and
  ``almost every'' in the topological sense: a condition is fulfilled
  almost everywhere (or in almost every point) if and only if the
  set of the points where it is fulfilled is dense in $M$.
\end{con}

\begin{proof}
  If $\lambda_i\ne 0$ at a point, then the Killing vector field $\bar
  \lambda^i $ is not identically zero. It is known that a Killing
  vector field that is not identically zero does not vanish on an
  open nonempty subset (to see it one can use the fact that the flow
  of a Killing vector field commutes with the exponential
  mapping). Thus, $\bar \lambda^i \ne 0$ at almost every point,
  implying $\lambda_i\ne 0$ at almost every point.
\end{proof}

\begin{cor}
  \label{cor:zero}
  Let $(M^{2n},g,J)$ be a connected Kähler manifold of dimension $2n\geq 4$ and let $(a_{ij},\lambda_i)$ be a solution of~\eqref{eq:hpr} such that
  $a_{i j}=0$ at every point of some open subset $U\subseteq M$. Then
  $(a_{ij},\lambda_i)\equiv(0,0)$ on the whole  $M$.
\end{cor}

\begin{proof}
  If $a_{ij}\equiv 0$ in $U$, then   $\lambda_i\equiv 0$ in  $U$ implying
  $\lambda_i\equiv 0$ on the whole $M$ in view of Corollary \ref{cor:killing2}. Then,  equation~\eqref{eq:hpr}
  implies that $a_{ij}$ is covariantly constant on $M$. Since it
  vanishes at a point, it vanishes everywhere. \end{proof}

\subsection{Algebraic lemma}
\par
Let us denote by $\mathcal J$ the following $(2,2)$-tensor:
\begin{equation}
\mathcal J^{\alpha \beta}_{i j} = \delta_i^\alpha \delta_j^\beta
+\J{\alpha}{i}{\beta}{j}.\label{eq:jten}
\end{equation}
Using this notation one can rewrite equation~(\ref{eq:hpr}) in the form
\begin{gather}
 \label{eq:hprJ}
 a_{i j,k}=\Jij(\lambda_\i g_{\j k} + \lambda_\j g_{\i k})
\end{gather}
The first step in the proof of Theorem \ref{thm:ext} will be to show the validity of the second equation of the system (\ref{eq:mg+Ba}) in a point:
\begin{lem}
  \label{lem:ext}
Let $(M^{2n},g,J)$ be a Kähler manifold of dimension $2n\geq 4$ and let $(a_{ij}, \lambda_i)$ and $(A_{ij}, \Lambda_i)$ be solutions of
  \eqref{eq:hpr} such that at the point $p\in M$, $a, g, $
  and $A$ are linearly independent.  Then, there exist numbers $B$ and
  $\mu$, such that the equation
  \begin{equation}
    \label{eq:mg+Ba1}
    \lambda_{i,j}=\mu g_{i j}+Ba_{i j}.
  \end{equation}
holds at $p$.
\end{lem}

\begin{proof}
  Substituting~(\ref{eq:hprJ}) in $a_{ij,kl} -
  a_{ij,lk}=a_{i \alpha}R^\alpha_{j k l}+a_{j \alpha}R^\alpha_{i k
    l},$ we obtain
  \begin{gather}
    \label{eq:int}
    a_{i \alpha}R^\alpha_{j k l}+a_{j \alpha}R^\alpha_{i k l}=\Jij
    (\lambda_{l, \i }g_{ \j k}+\lambda_{l, \j}g_{ \i k}-\lambda_{k, \i }g_{ \j l}-\lambda_{k, \j}g_{ \i l})
  \end{gather}
  These equations are fulfilled for every solution of~(\ref{eq:hpr}),
  thus for $(A_{ij}, \Lambda_i)$. We denote by~(\ref{eq:int}.A) the
  equation~(\ref{eq:int}) with $(a_{i j},\lambda_i)$ replaced by
  $(A_{i j},\Lambda_i)$.  From this point we will 
  work in the tangent space to the  fixed point $p$ only.

  \par
  Since equations~(\ref{eq:int}) and~(\ref{eq:int}.A) are not affected
  by the transformation (for any constants $a, A, c, C$)
  \begin{gather}
    a_{i j}\to a_{i j}+a\cdot g_{i j},\qquad \lambda_{i,j}\to\lambda_{i,j}+c\cdot g_{ij},\\
    A_{i j}\to A_{i j}+A\cdot g_{i j},\qquad \Lambda_{i,j}\to\Lambda_{i,j}+C\cdot g_{ij},
  \end{gather}
  without loss of generality we can assume that $a_{i j}$,
  $\lambda_{i, j}$, $A_{i j}$ and $\Lambda_{i,j}$ are trace-free, i.e.
  \begin{gather}
    a_{ij}g^{ij}=\lambda_{i,j}g^{ij}=A_{ij}g^{ij}=\Lambda_{i,j}g^{ij}=0.
  \end{gather}
    In this ``trace-free'' situation, our goal is to show that
  $\lambda_{i,j} = B \cdot a_{ij}$ for a certain number $B$.

  \par
  After contracting~(\ref{eq:int}) with $A_\l^l$ and renaming of the indices
  $l\to\beta$, $\l\to l$, we obtain:
  \begin{multline}
    \label{eq:Riem}
    a_{i \alpha} R^\alpha_{j k \beta} A_l^\beta + a_{j \alpha}
    R^\alpha_{i k \beta} A_l^\beta
    =\Jij (A_l^\beta \lambda_{\beta, \i }g_{ \j k}+ A_l^\beta
    \lambda_{\beta, \j}g_{ \i k}- A_l^\beta \lambda_{k, \i}g_{\j
      \beta}- A_l^\beta \lambda_{k, \j}g_{ \i \beta}).
  \end{multline}
  Because of the symmetries of the curvature tensor,
  $$a_{i \alpha} R^\alpha_{j k \beta} A_l^\beta=a_i^\alpha R_{\alpha j k \beta} A_l^\beta=a_i^\alpha R_{\beta k j \alpha } A_l^\beta.$$
  Then, equation~(\ref{eq:Riem}) can be rewritten as
  \begin{multline}
    a^\alpha_i A_l^\beta R_{\beta k j \alpha}+a^\alpha_j A_l^\beta
    R_{\beta k i \alpha}
    =\Jij (A_l^\beta \lambda_{\beta, \i }g_{ \j k}+A_l^\beta
    \lambda_{\beta, \j}g_{ \i k}-A_l^\beta \lambda_{k, \i}g_{ \j
      \beta}-A_l^\beta \lambda_{k, \j}g_{ \i \beta}).
  \end{multline}
  Symmetrizing with respect to  $(l,k)$ and rearranging the terms we obtain
  \begin{multline}
    \label{11}
    a^\alpha_i(A_l^\beta R_{\beta k j \alpha}+ A_k^\beta R_{\beta l j
      \alpha})+a^\alpha_j (A_l^\beta
    R_{\beta k i \alpha}+ A_k^\beta R_{\beta l i \alpha})=\\
    =\Jij (A_l^\beta \lambda_{\beta, \i }g_{ \j k}+A_l^\beta
    \lambda_{\beta, \j}g_{ \i k} -A_l^\beta \lambda_{k, \i }g_{ \j
      \beta}-A_l^\beta \lambda_{k, \j}g_{ \i \beta} + \\
    + A_k^\beta \lambda_{\beta, \i }g_{ \j l}+A_k^\beta
    \lambda_{\beta, \j}g_{ \i l}-A_k^\beta \lambda_{l, \i }g_{ \j
      \beta}-A_k^\beta \lambda_{l, \j}g_{ \i \beta}).
  \end{multline}
  The terms in the brackets in the left hand side are the left hand
  side of~(\ref{eq:int}.A) with renamed indices: the rules for
  renaming indices are
  $$
  \begin{pmatrix}
    i & \alpha & j&k&l\\
    l & \beta & k&j&\alpha
  \end{pmatrix}
  \mbox{ and }
  \begin{pmatrix}
    i & \alpha & j&k&l\\
    l & \beta & k&i&\alpha
  \end{pmatrix},
  $$
  respectively. Substituting~(\ref{eq:int}.A) in \eqref{11}, we obtain
  \begin{multline}
    \label{eq:Jkl}
    \Jkl \left[ a^\alpha_i (\Lambda_{ \alpha, \l }g_{ \k j}+\Lambda_{ \alpha, \k}g_{ \l j}-\Lambda_{j, \l }g_{ \k \alpha }-\Lambda_{j, \k}g_{ \l \alpha })+\right.\\
    \left.
      +a^\alpha_j
      (\Lambda_{ \alpha, \l }g_{ \k i}+\Lambda_{ \alpha, \k}g_{ \l i}-\Lambda_{i, \l }g_{ \k \alpha }-\Lambda_{i, \k}g_{ \l \alpha })\right] ={}\\
    =\Jij \left[ A_l^\alpha (\lambda_{\alpha, \i }g_{ \j k}+\lambda_{\alpha, \j}g_{ \i k}-\lambda_{k, \i }g_{ \j \alpha}- \lambda_{k, \j}g_{ \i \alpha})+\right.\\
    \left. +A_k^\alpha (\lambda_{\alpha, \i }g_{ \j l}+\lambda_{\alpha, \j}g_{ \i l}-\lambda_{l, \i }g_{ \j \alpha}-\lambda_{l, \j}g_{ \i \alpha})\right].
  \end{multline}
  Now we want to change the contraction with the tensor $\Jkl$ by the
  contraction with the tensor $\Jij$.  This operation is possible (=
  after applying it we obtain the same  equation), because of
  specific symmetries of each component in brackets. Indeed, for the
  first component we have
  \begin{multline}
    \Jkl a^\alpha_i \Lambda_{ \alpha, \l}g_{ \k j}=(\delta_k^\k\delta_l^\l+\J{\k}{k}{\l}{l}) a^\alpha_i
    \Lambda_{ \alpha, \l}g_{ \k j}=\\
    =\delta_k^\k\delta_l^\l a^\alpha_i \Lambda_{ \alpha, \l}g_{ \k j}+\J{\k}{k}{\l}{l} a^\alpha_i \Lambda_{ \alpha, \l}g_{ \k j}=\\
    =\delta_i^\i\delta_j^\j a^\alpha_\i \Lambda_{ \alpha, l}g_{ k \j}+\J{\k}{k}{\l}{l} a^\alpha_i \Lambda_{ \alpha, \l}g_{ \k j}
  \end{multline}
  Consider the last part and apply Lemma~\ref{lem:skew} several times:
  \begin{multline}
    \J{\k}{k}{\l}{l} a^\alpha_i \Lambda_{ \alpha, \l}g_{ \k j}=
    (J_{\ \ k}^\k g_{ \k j})\cdot (J_{\ \ l}^\l \Lambda_{\l,\alpha})\cdot(g^{\alpha \beta}a_{\beta i})=\\
    =(-J_{\ \ j}^\j g_{ \j k})\cdot (-J_{\ \ \alpha}^{\alpha'} \Lambda_{\alpha',l})\cdot (g^{\alpha\beta} a_{\beta i})=
    (J_{\ \ j}^\j g_{ \j k})\cdot\Lambda_{\alpha', l}\cdot(
    J_{\ \ \alpha}^{\alpha'} g^{\alpha \beta})\cdot a_{\beta i}=\\
    =(J_{\ \ j}^\j g_{ \j k})\cdot\Lambda_{\alpha',l}\cdot(-J_{\ \ \beta'}^{\beta}g^{\beta' \alpha'})\cdot a_{\beta i}=
    (J_{\ \ j}^\j g_{ \j k})\cdot\Lambda_{\alpha', l}g^{\beta'\alpha'}\cdot(-J_{\ \ \beta'}^{\beta} a_{\beta i})=\\
    =(J_{\ \ j}^\j g_{ \j k})\cdot\Lambda_{\alpha', l}g^{\beta'\alpha'}\cdot(J_{\ \ i}^{\i} a_{\i \beta'})
    =J_{\ \ i}^{\i} J_{\ \ j}^\j a_{\i}^{\alpha} \Lambda_{\alpha, l} g_{ \j k}
  \end{multline}
  Then $\Jkl a^\alpha_i \Lambda_{ \alpha, \l}g_{ \k j} = \Jij
  a^\alpha_\i \Lambda_{ \alpha, l}g_{ k \j}$, as we claimed.

  \par
  The proof for all other components is analogous (in fact, in the
  proof we used the hermitian property of $a_{ij}, \lambda_{i,j},
  g_{ij}$ only, and this property is fulfilled for all these tensors
  by Lemma~\ref{lem:skew}).

  \par
  Therefore, considering each component in the left part
  of~(\ref{eq:Jkl}) separately, we obtain:
  \begin{multline}
    \label{eq:Jij}
    \Jij \left[ a^\alpha_\i (\Lambda_{ \alpha, l }g_{ k \j}+\Lambda_{ \alpha, k}g_{ l \j}-\Lambda_{\j, l }g_{k \alpha }-\Lambda_{\j, k}g_{ l \alpha })+\right.\\
    \left.
      +a^\alpha_\j
      (\Lambda_{ \alpha, l }g_{ k \i}+\Lambda_{ \alpha, k}g_{ l \i}-\Lambda_{\i, l }g_{ k \alpha }-\Lambda_{\i, k}g_{ l \alpha })\right] ={}\\
    =\Jij \left[ A_l^\alpha (\lambda_{\alpha, \i }g_{ \j k}+\lambda_{\alpha, \j}g_{ \i k}-\lambda_{k, \i}g_{ \j \alpha}- \lambda_{k, \j}g_{ \i \alpha})+\right.\\
    \left. +A_k^\alpha (\lambda_{\alpha, \i }g_{ \j l}+\lambda_{\alpha, \j}g_{ \i l}-\lambda_{l, \i }g_{ \j \alpha}-\lambda_{l, \j}g_{ \i \alpha})\right].
  \end{multline}
  In the left hand side of \eqref{eq:Jij}, we  collect the components on containing $g$ with the same indices: 
  \begin{multline}
    \label{16}
    \Jij\left[
      (a^\alpha_\i\Lambda_{\alpha, l}-A_l^\alpha\lambda_{\alpha, \i})g_{k \j}
      +
      (a^\alpha_\i\Lambda_{\alpha, k}-A_k^\alpha\lambda_{\alpha, \i})g_{l \j}
      +
    \right.\\\left.
      +
      (a^\alpha_\j\Lambda_{\alpha, l}-A_l^\alpha\lambda_{\alpha, \j})g_{k \i}
      +
      (a^\alpha_\j\Lambda_{\alpha, k}-A_k^\alpha\lambda_{\alpha, \j})g_{l \i}
    \right]={}\\
    =\Jij\left[
      a_{\i k}\Lambda_{\j, l}+a_{\i l}\Lambda_{\j, k}
      +a_{\j k}\Lambda_{\i, l}+a_{\j l}\Lambda_{\i, k}
      -A_{\j l}\lambda_{k, \i}-A_{\i l}\lambda_{k, \j}
      -A_{\j k}\lambda_{l, \i}-A_{\i k}\lambda_{l, \j}
    \right]
  \end{multline}
  We set  
  $c_{i l}=a^\alpha_i \Lambda_{\alpha, l}-A^\alpha_l\lambda_{\alpha, i}.$
  it is easy to check that $c_{il}$ anticommutes with $J$:
  $ J_{\ \ i}^\alpha c_{\alpha j}=-c_{i \alpha} J_{\ \ j}^\alpha.$
  Then  equation~\eqref{16} takes the form:
  \begin{multline}
    \label{eq:c}
    \Jij\left[c_{\i l}g_{\j k}+c_{\i k}g_{\j l}+c_{\j l}g_{\i k}+c_{\j k}g_{\i l}
    \right]=\\
    =\Jij\left[
      a_{\i k}\Lambda_{\j, l}+a_{\i l}\Lambda_{\j, k}
      +a_{\j k}\Lambda_{\i, l}+a_{\j l}\Lambda_{\i, k}
      -A_{\j l}\lambda_{k, \i}-A_{\i l}\lambda_{k, \j}
      -A_{\j k}\lambda_{l, \i}-A_{\i k}\lambda_{l, \j}
    \right]
  \end{multline}
  Let us now contract the last equation with $g^{j k}$. This operation
  involves the $j$-index, so we have to make use of the explicit formula (\ref{eq:jten})
  for $\mathcal J$. After some index manipulations, using the
  anticommutation- and trace-free-properties of the tensors
  involved\footnote{ Each component separately:
    \begin{gather*}
      c_{i l}g_{jk} g^{jk}= 2nc_{il},
      \quad
      c_{i k}g_{jl} g^{jk}= c_{il},\\
      c_{j l}g_{ik} g^{jk}= c_{il},
      \quad
      c_{j k}g_{il} g^{jk}= (c_{jk}g^{jk})g_{il},\\
      \J{\i}{i}{\j}{j}c_{\i l}g_{\j k} g^{jk}= 0
      % \mbox{ (trace of skew-symmetric form)}
      ,
      \quad
      \J{\i}{i}{\j}{j}c_{\i k}g_{\j l} g^{jk}= -c_{il},\\
      \J{\i}{i}{\j}{j}c_{\j l}g_{\i k} g^{jk}= -\J{\i}{i}{\l}{l}c_{\i \l}=-c_{il},
      \quad
      \J{\i}{i}{\j}{j}c_{\j k}g_{\i l} g^{jk}= -(g^{jk}J^{\j}_{\ \ j}c_{\j k})J^{\i}_{\ \ i}g_{\i l}=0,\\
      a_{i k}\Lambda_{j, l} g^{jk}= a_i^p\Lambda_{p,l},
      \quad
      a_{i l}\Lambda_{j,k} g^{jk}= 0,
      \quad
      a_{j k}\Lambda_{i, l} g^{jk}= 0,
      \quad
      a_{jl}\Lambda_{i, k} g^{jk}= \Lambda_{i, p} a^p_l,\\
      \J{\i}{i}{\j}{j}a_{\i k}\Lambda_{\j, l} g^{jk}= -a_i^p\Lambda_{p,l},
      \quad
      \J{\i}{i}{\j}{j}a_{\i l}\Lambda_{\j, k} g^{jk}= 0,
      \quad
      \J{\i}{i}{\j}{j}a_{\j k}\Lambda_{\i, l} g^{jk}= 0 ,
      \quad
      \J{\i}{i}{\j}{j}a_{\j l}\Lambda_{\i, k} g^{jk}= -\Lambda_{i, p} a^p_l.\\
    \end{gather*}}
, we obtain:
  \begin{gather}
    \label{eq:cil}
    2n c_{il}+(c_{jk}g^{jk})g_{il}=0,
  \end{gather}
  which implies $c_{il}=0$. Since $c_{il}=0$, the
  equation~(\ref{eq:c}) reads
  \begin{equation}
    \label{eq:withoutc}
    \Jij\left[
      a_{\i k}\Lambda_{\j, l}+a_{\i l}\Lambda_{\j, k}
      +a_{\j k}\Lambda_{\i, l}+a_{\j l}\Lambda_{\i, k}
      -A_{\j l}\lambda_{k, \i}-A_{\i l}\lambda_{k, \j}
      -A_{\j k}\lambda_{l, \i}-A_{\i k}\lambda_{l, \j}
    \right]=0
  \end{equation}
  Let us now multiply~(\ref{eq:withoutc}) by $\frac{1}{2}\mathcal J_{p
    q}^{j k}$. After rearranging components and renaming indices we
  can write the equation in a more symmetric way:
  \begin{multline}
    \frac{1}{2}(\delta_i^\i \delta_j^\j \delta_k^\k+\delta_i^\i J_{\ \ j}^\j J_{\ \ k}^\k+
    J_{\ \ i}^\i J_{\ \ j}^\j \delta_k^\k - J_{\ \ i}^\i \delta_j^\j J_{\ \ k}^\k)\cdot\\
    \cdot(
    a_{\i \k}\Lambda_{\j, l}+a_{\i l}\Lambda_{\j, \k}
    +a_{\j \k}\Lambda_{\i, l}+a_{\j l}\Lambda_{\i, \k}-\\
    -A_{\j l}\lambda_{\k, \i}-A_{\i l}\lambda_{\k, \j}
    -A_{\j \k}\lambda_{l, \i}-A_{\i \k}\lambda_{l, \j}
    )=0
  \end{multline}
  Using that $J$ anticommutes with $a_{ij}$, $A_{ij}$, $\lambda_{i,j}$
  (see Lemma~\ref{lem:skew}) one can get
  \begin{gather}
    \label{eq:Jsym}
    a_{i l}\Lambda_{j, k}
    +a_{j k}\Lambda_{i, l}+
    \J{\i}{i}{\j}{j}(a_{\i l}\Lambda_{\j, k}
    +a_{\j k}\Lambda_{\i, l})
    =A_{i l}\lambda_{j, k}
    +A_{j k}\lambda_{i, l}+
    \J{\i}{i}{\j}{j}(A_{\i l}\lambda_{\j, k}
    +A_{\j k}\lambda_{\i, l})
  \end{gather}
  Symmetrizing~(\ref{eq:Jsym}) by $(i, l)$ we finally obtain
  \begin{gather}
    a_{i l}\Lambda_{j, k}
    +a_{j k}\Lambda_{i, l}
    =A_{i l}\lambda_{j, k}
    +A_{j k}\lambda_{i, l}.
  \end{gather}
  In other words, $\Lambda_{\aleph} a_{\daleth} + \Lambda_{\daleth}
  a_{\aleph} = \lambda_{ \daleth} A_{\aleph} +\lambda_{ \aleph}
  A_{\daleth}$, where $\aleph$ and $\daleth$ stand for the symmetric
  indices $jl$ and $ik$, respectively.

  \par
  But it is easy to check that a non-zero simple symmetric tensor
  $X_{\aleph \daleth} = P_\aleph Q_\daleth + P_\daleth Q_\aleph$ determines
  its factors $P_\aleph$ and $Q_\daleth$ up to scale and order (it is
  sufficient to check, for example, by taking $P_\aleph$ and $Q_\daleth$
  to be basis vectors). Since $a_{ij}$ and $A_{ij}$ are supposed to be
  linearly independent, it follows that $\lambda_{i,j} = \const \cdot
  a_{ij}$, as required.
\end{proof}

\begin{rem}
  We would like to emphasize here that, though Lemma~\ref{lem:ext} is
  formulated in the differential-geometrical notation, it is
  essentially an algebraic statement (in the proof we did not use
  differentiation except for the integrability conditions \eqref{eq:int}
  that were actually obtained before, see \eqref{integr}). Moreover,
  we can replace $R^i_{j kl}$ in \eqref{eq:int} by any $(1,3)$-tensor
  having the same algebraic symmetries (with respect to $g$) as the
  curvature tensor.
\end{rem}

\subsection{ If the solutions $a_{ij}, A_{ij}$ and $g_{ij}$ are
  linearly dependent over functions, then they are linearly dependent
  over constants}
\label{sec:lindep}
The goal of this section is to show, that under the assumption of degree of mobility $\geq 3$, equation (\ref{eq:mg+Ba1}) holds in a neighborhood of almost every point of $M$ for each solution $(a_{ij},\lambda_{i})$ of equation (\ref{eq:hpr}). The real numbers $B$ and $\mu$ in equation (\ref{eq:mg+Ba1}) then become smooth function on this neighborhood. In the end of this section, it will be also shown that the local function $B$ is the same for all solutions of equation (\ref{eq:hpr}). 
\begin{lem}
  \label{cor:lindep}
On a Kähler manifold $(M^{2n\ge 4} ,g,J)$, let $(A_{ij}, \lambda_i)$ and $(a_{ij}, \lambda_i)$ be solutions
  of~(\ref{eq:hpr}). Then, almost every point $p\in M$ has a
  neighborhood $U(p)\ni p$ such that in this neighborhood one of the
  following conditions is fulfilled:

  \begin{itemize}
  \item[(a)] $a_{ij}, A_{ij}, $ and $g_{ij}$ are linearly independent
    at every point of $U(p)$,
  \item[(b)] $a_{ij}, A_{ij}, $ and $g_{ij}$ are linearly dependent at
    every point of $U(p)$.
  \end{itemize}
\end{lem}

\begin{proof} The proof in fact 
does not require that $a_{ij}$ and  $A_{ij}$ are solutions of (\ref{eq:hpr}). 
  Let $W$ be the set of the points where (a) is fulfilled. $ W$ is
  evidently an open set. Consider $\mbox{int} \left( M\setminus
    W\right)$, where ``$\mbox{int}$'' denotes the set of the interior
  points. This is also an open set, and $W \cup \mbox{int}\left( M\setminus
    W\right) $ is open and everywhere dense. By construction, every
  point of $W \cup \mbox{int}\left( M\setminus W\right) $ has a neighborhood
  satisfying the condition (a) or the condition (b). \end{proof}

One of the possibilities in Lemma~\ref{cor:lindep} is that (in a
neighborhood $U(p)$) the solutions $a_{ij}, A_{ij}$ and $g_{ij}$ of
\eqref{eq:hpr} are linearly depended over functions. Our goal is to
show that in this case they are actually linearly dependent (over
constants).  At first we consider the special  case, 
when  two solutions are  proportional.

\begin{lem}
  \label{lem:A=ca}
  Let $(M^{2n},g,J)$ be a Kähler manifold  of dimension $2n\ge 4$,  
  and let $(a_{ij},\lambda_{i})$ and $(A_{ij},\Lambda_{i})$ be solutions
  of~(\ref{eq:hpr}) such that $a_{i j}\neq 0$ at every point of some open subset $U\subseteq M$.
  If $\alpha:U\to \R$ is a function such that
  \begin{gather}
    \label{eq:A=ca}
    A=\alpha\,a, 
  \end{gather}
  then $\alpha$ is  constant, and $A=\alpha\,a$ on the whole $M$.
\end{lem}

\begin{proof}
  Since $A_{ij}$ and $a_{ij}$ are smooth tensor fields on $U$ and
  $a_{i j}\neq 0$, the function $\alpha$ is also smooth.  We
  covariantly differentiate~\eqref{eq:A=ca} and substitute the
  derivatives of $a_{ij}$ and $A_{ij}$ using~(\ref{eq:hpr}) to obtain
  \begin{align}
    \gamma_{i}g_{jk}+\gamma_{j}g_{ik}-\bar{\gamma}_{i}J_{jk}-\bar{\gamma}_{j}J_{ik}=\alpha_{,k}a_{ij},
    \label{eq:globind1}
  \end{align}
  where $\gamma_{i}:=\Lambda_{i}-\alpha \lambda_{i}$. Contracting
  equation~(\ref{eq:globind1}) with a non-zero vector field $U^{k}$
  such that $U^k \alpha_{,k}=0$ yields
  \begin{equation}
    \gamma_{i}U_{j}+\gamma_{j}U_{i}+\bar{\gamma}_{i}\bar{U}_{j}+\bar{\gamma}_{j}\bar{U}_{i}=0
    \label{eq:globind2}
  \end{equation}
  Let us now show that at every point
  $$\mbox{span}\{U^{j},\bar{U}^{j}\}^{\perp}\subseteq \mbox{span}\{\gamma^{j},\bar{\gamma}^{j}\}^{\perp}.$$
  For every vector field
  $V^{j}\in\mbox{span}\{U^{j},\bar{U}^{j}\}^{\perp}$ we have
  (contracting this vector field with \eqref{eq:globind2})
  \begin{align}
    (\gamma_{j}V^{j})\, U_{i}+ (\bar{\gamma}_{j}V^{j})\, \bar{U}_{i}=0\nonumber
  \end{align}
  Since $U_i$ and $\bar U_i$ are linearly independent,
  $\gamma_{j}V^{j}=\bar\gamma_{j}V^{j}=0$. Then $V_i \in
  \mbox{span}\{\gamma^{j},\bar{\gamma}^{j}\}^\perp$. Thus,
  $\mbox{span}\{U^{j},\bar{U}^{j}\}^\perp\subseteq\mbox{span}\{\gamma^{j},\bar{\gamma}^{j}\}^\perp$ as we claimed. 

  \par
  Assume $\gamma_i\neq 0$. Then the spaces
  $\mbox{span}\{U^{j},\bar{U}^{j}\}^{\perp}$ and
  $\mbox{span}\{\gamma^{j},\bar{\gamma}^{j}\}^{\perp}$ have equal
  dimension $(2n-2)$, and therefore coincide. The same holds for their orthogonal
  complements and we obtain
  $$\mbox{span}\{U^{j},\bar{U}^{j}\}=\mbox{span}\{\gamma^{j},\bar{\gamma}^{j}\}$$
  Thus, every vector $U^i$ from the at least $(2n-1)$-dimensional space
  $\mbox{span}(\alpha_{,}^{\ \ i})^\perp$ lies in the $2$-dimensional
  space $\mbox{span}\{\gamma^{j},\bar{\gamma}^{j}\}$, which gives us a
  contradiction. Thus, $\gamma_i=0$ and equation~(\ref{eq:globind1})
  reads $\alpha_{,k} a_{ij}=0$, implying $\alpha$ is constant on $U$.   Therefore, the solution $A_{ij}- \alpha  a_{ij}$ vanishes at
  every point of $U$. By Corollary~\ref{cor:zero} it vanishes on the
  whole  $M$.
\end{proof}
Now let us treat the general case:
\begin{lem}
  \label{lem:globind}
  On a connected Kähler manifold $(M^{2n},g,J)$ of dimension $2n\geq 4$, let $(a_{ij},\lambda_{i})$ and $(A_{ij},\Lambda_{i})$ be solutions
  of~(\ref{eq:hpr}). Assume that for certain functions $\alpha$ and $\beta$ on
  an open subset $U\subseteq M$ we have
  \begin{gather}
    \label{eq:A=g+a}
    A_{ij}=\alpha g_{ij}+\beta a_{ij}
  \end{gather}
  Then there exist constants $(C_1, C_2, C_3)\neq(0,0,0)$ such that
  $$C_1 A+C_2 a+ C_3 g=0 \quad \mbox{on the whole  $M$.}$$
\end{lem}

\begin{proof}
  If there locally exists a function $c$ such that $a_{i j}=c g_{i j}$, then by
  the previous Lemma~\ref{lem:A=ca} the function $c$ is a constant. Hence,  by Corollary~\ref{cor:zero}, one can choose
  $C_1=0$, $C_2=-1$ and $C_3=c$.

  \par Let $a_{i j}$ be non-proportional to
  $g_{ij}$. Then~\eqref{eq:A=g+a} is a linear system of equations of
  maximal rank with smooth coefficients on functions $\alpha$ and
  $\beta$. Thus, its solutions $\alpha$ and $\beta$ are smooth.

  \par Similarly as before in Lemma~\ref{lem:A=ca}, by differentiating \eqref{eq:A=g+a}
  we obtain
  \begin{gather}
    \label{eq:globind3}
    \gamma_{i}g_{jk}+\gamma_{j}g_{ik}-\bar{\gamma}_{i}J_{jk}-\bar{\gamma}_{j}J_{ik}=\alpha_{,k}g_{ij}+\beta_{,k}a_{ij}
  \end{gather}
  where $\gamma_i = \Lambda_{i}-\beta \lambda_{i}$.

  \par Assume $\gamma_i\neq 0$. We contract \eqref{eq:globind3} with a vector field
  $U^i$ such that  $U^k \alpha_{,k}= U^k\beta_{,k}=0$ to obtain
  equation~\eqref{eq:globind2}. As in  the proof of
  Lemma~\ref{lem:A=ca},  we obtain
  $$\mbox{span}\{U^{j},\bar{U}^{j}\}=\mbox{span}\{\gamma^{j},\bar{\gamma}^{j}\}$$
  implying $U_i=c\cdot\gamma_i+d\cdot\bar\gamma_i$ for certain
  functions $c$ and $d$. We substitute $U_i$ in~\eqref{eq:globind2} to
  obtain
  $$2 c\cdot(\gamma_i \gamma_j+\bar\gamma_i \bar\gamma_j)=0.$$
  Since $\gamma_i\neq0$, it follows that $c=0$, and therefore $U^i=d\cdot\bar\gamma^i$. We have shown that every vector $U^i$ from
  the at least $(2n-2)$-dimensional space
  $\mbox{span}(\alpha_{,}^{\ i},\beta_{,}^{\ i})^\perp$ is proportional to $\bar
  \gamma^i$, which gives us a contradiction. Thus, $\gamma_i=0$ and equation (\ref{eq:globind3}) takes the form
  $$\alpha_{,k}a_{ij}+\beta_{,k}g_{ij}=0.$$
  We have $\alpha_{,k}\equiv 0\equiv \beta_{,k}$, implying
  $\alpha\equiv \const=:C_{2}$ and $\beta=\const=:C_{3}$.

  \par Therefore, the solution $A_{ij}- C_2 a_{ij}-C_3 g_{ij}$ vanishes at
  every point of $U$. By Corollary~\ref{cor:zero} it vanishes on the
  whole  $M$.
\end{proof}

Thus, if the degree of mobility is $\geq 3$, by Lemma~\ref{lem:globind},
for every  solution $(a_{ij},\lambda_{i})$ of 
(\ref{eq:hpr})  such that $a_{ij}\neq \mbox{const}\cdot g_{ij}$, 
 equation~(\ref{eq:mg+Ba1}) holds in a neighborhood of almost every
point of $M$ (for some locally defined functions $B$ and $\mu$ that could a priori 
depend on the solution $(a_{ij},\lambda_{i})$). Our next goal is to 
show, that the function $B$ is the same for all  solutions:

\begin{cor}
  \label{cor:Bind}
  Let $(M^{2n},g,J)$ be a Kähler manifold of dimension $2n\geq 4$ and assume that the degree of mobility is $\geq 3$. Then, the function $B$ defined by
  equation~(\ref{eq:mg+Ba1}) does not depend on the solution
  $(a_{ij},\lambda_{i})$ of equation~(\ref{eq:hpr}).
\end{cor}

\begin{proof}
  Take the  second solution $(A_{ij},\Lambda_{i})$ of
  equation~(\ref{eq:hpr}).  Let us  first  assume that $g_{ij}$, $a_{ij}$
  and $A_{ij}$ are linearly independent.

  \par We know that $(a_{ij}+A_{ij},\lambda_{i}+\Lambda_{i})$ is again
  a solution. Adding  equations~(\ref{eq:mg+Ba1}) for
  $(a_{ij},\lambda_{i})$ and $(A_{ij},\Lambda_{i})$ with functions $B$
  and $B'$ respectively and substracting the same equation
  corresponding to the sum of the both solutions (the correspondent  function $B$ for the sum of solutions  will be  denoted by $B^{+}$), we obtain 
  \begin{align}
    0=\mbox{something}\cdot g_{ij}+(B-B^{+})a_{ij}+(B'-B^{+})A_{ij}\nonumber
  \end{align}
 Combining  Lemma~\ref{lem:globind} and the assumption that $g, a$ and $A$ are linearly independent, we obtain   $B= B^+=B'$  as we claimed.

  \par Consider now the second  case when $g_{ij}$, $a_{ij}$ and $A_{ij}$ are
  linearly dependent, i.e. (without loss of generality), 
   $A_{ij}=C g_{ij}+D
  a_{ij}$ on $M$ for some constants $C$ and $D$. Thus, the
  corresponding $1$-forms $\Lambda_{i}$ and $\lambda_{i}$ for $A_{ij}$
  and $a_{ij}$ respectively are related by the equation $\Lambda_{i}=D
  \lambda_{i}$. Multiplying   equation (\ref{eq:mg+Ba1}) by $D$ we obtain 
  \begin{gather}
    \underbrace{D \lambda_{i,j}}_{\Lambda_{i,j}}=D\mu g_{ij} + D B
    a_{ij} =\underbrace{(D \mu - C B)}_{\Mu} g_{ij} + \underbrace{(D
      a_{ij} + C g_{i j})}_{A_{ij}} B
  \end{gather}
  This is   equation (\ref{eq:mg+Ba1}) on $(A_{ij}, \Lambda_i)$ with the same
  function $B$. Finally, in all cases,  the function $B$ is the same for all solutions of equation (\ref{eq:hpr}).
\end{proof}

\subsection{In the neighborhood of a point such that $g, a, $ and
  $A$ are linearly independent, the function $B$ is a constant}

Our next goal is to show that the local function $B$ we have found is a constant.

\begin{lem}
  \label{lem:ext2}
Let $(M^{2n},g,J)$ be a Kähler manifold of dimension $2n\geq 4$.  
Suppose that in a neighborhood $U\subseteq M$ there exist at least
  two solutions $(a_{ij}, \lambda_i)$ and $(A_{ij}, \Lambda_i)$
  of~(\ref{eq:hpr}) such that $a, A$ and $g$ are linearly independent
  at every point of $U$.  Then the function $B$ defined by
  equation~(\ref{eq:mg+Ba1}) is a constant.
\end{lem}

The proofs for the cases  $\dim M\geq 6$ and $\dim M=4$ use different
methods and will be given in sections~\ref{Bconstdim6}
and~\ref{Bconstdim4} respectively.
\subsubsection{Proof of Lemma~\ref{lem:ext2}, if  $\dim M\geq 6$.}
\label{Bconstdim6}
First of all, the function $B$ is smooth. Indeed, the trace-free
version of \eqref{eq:mg+Ba1} is
\begin{equation}
  \label{tracefree1}
  \lambda_{i,j}-\tfrac{1}{2n}\lambda_{k,}^{\ \ k} \cdot g_{ij}= B(a_{ij} - \tfrac{2}{n} \lambda g_{ij}),
\end{equation}
where $\lambda:= \tfrac{1}{4}a^i_i$, and the function $B$ is smooth since it is the coefficient of the
proportionality of the nowhere vanishing tensor $(a_{ij} -
\tfrac{2}{n} \lambda g_{ij})$ and the tensor $(
\lambda_{i,j}-\tfrac{1}{2n}\lambda_{k,}^{\ \ k}\cdot  g_{ij})$.  Since $B$ is
smooth, $\mu$ is smooth as well, as the coefficient of the
proportionality of the nowhere vanishing tensor $g_{ij}$ and the
tensor $( \lambda_{i,j}- Ba_{i j})$.

\par Thus, all objects in the equation
\begin{gather}
  \label{eq:lmuB}
  \lambda_{i,j}=\mu g_{i j}+Ba_{i j}
\end{gather}
are smooth. We covariantly differentiate the equation and substitute
$a_{ij,k}$ using~(\ref{eq:hprJ}) to obtain
\begin{equation} \label{-2}
  \lambda_{i,jk}=\mu_{,k}g_{i j} + B_{,k}a_{i j}+B a_{ij,k}
  =\mu_{,k}g_{i j} + B_{,k}a_{i j}+B\cdot \Jij(\lambda_\i g_{\j k} +
  \lambda_\j g_{\i k}).
\end{equation}

By definition of the curvature tensor,
\begin{multline}
  \label{eq:lR}
  \lambda_p R^p_{ijk}=\lambda_{i,jk}-\lambda_{i,kj}
  \stackrel{\eqref{-2}}{=}\mu_{,k}g_{i j}-\mu_{,j}g_{i k}+B_{,k}a_{i j}- B_{,j}a_{i k}+\\
  +B\cdot \Jij(\lambda_\i g_{\j k} + \lambda_\j g_{\i k})-B\cdot
  \mathcal{J}_{i k}^{\i\k}(\lambda_\i g_{\k j} + \lambda_\k g_{\i j})=\\
  =\mu_{,k}g_{i j}-\mu_{,j}g_{i k}+B_{,k}a_{i j}- B_{,j}a_{i k}+
  B \lambda_j g_{i k}- B \lambda_k g_{i j}+\\
  +B\cdot \J{\i}{i}{\j}{j}(2 \lambda_\i g_{\j k} + \lambda_\j g_{\i k})- B \J{\i}{i}{\k}{k}\lambda_\k g_{\i j}
\end{multline}

Let us now substitute $\lambda_{i, j}$ in~(\ref{eq:int})
by~(\ref{eq:lmuB}). The components with $\mu$ disappear because of
the symmetries of $g_{ij}$ and the equation takes the following form:
\begin{gather}
  \label{eq:int_cond}
  a_{i \alpha}R^\alpha_{j k l}+a_{j \alpha}R^\alpha_{i k l}
  =B \Jij (a_{l \i }g_{ \j k}+a_{l \j}g_{ \i k}-a_{k \i }g_{ \j l}-a_{k \j}g_{ \i l})
\end{gather}

We contract this equation with $\lambda^l$. Applying the identity
$a_{i \alpha}R^\alpha_{jk\beta}\lambda^\beta=a_i^\alpha \lambda_\beta
R^\beta_{kj\alpha}$ we obtain
\begin{gather}
  a_i^\alpha \lambda_\beta R^\beta_{kj\alpha}+
  a_j^\alpha \lambda_\beta R^\beta_{ki\alpha}
  =B \Jij (\lambda^\beta a_{\beta \i }g_{ \j k}
  +\lambda^\beta a_{\beta \j}g_{ \i k}
  -a_{k \i }\lambda_\j -a_{k \j}\lambda_\i).
\end{gather}

Now we substitute the left hand side using~(\ref{eq:lR}). After substituting (\ref{eq:jten}) for $\Jij$ and tensor manipulation, we
 obtain
\begin{multline}
  g_{k j} (a_i^\alpha \mu_{,\alpha}-2 B\lambda^\alpha a_{i \alpha})+
  g_{k i} (a_j^\alpha \mu_{, \alpha}-2 B\lambda^\alpha a_{j \alpha})+\\
  +
  a_{k j} (a_i^\alpha B_{,\alpha}-\mu_{, i} + 2 B\lambda_i)+
  a_{k i} (a_j^\alpha B_{,\alpha}-\mu_{, j} + 2 B\lambda_j)=\\
  =B_{,j} a_{k \alpha}a^\alpha_i + B_{,i} a_{k\alpha}a^\alpha_j
\end{multline}

Set  $\xi_i:=a_i^\alpha \mu_{, \alpha}-2 B\lambda^\alpha a_{i \alpha}$
and $\eta_i:=a_i^\alpha B_{,\alpha}-\mu_{, i} + 2 B\lambda_i$.
Then
\begin{gather}
  \label{eq:xe}
  \xi_i g_{k j}+ \xi_j g_{k i}+
  \eta_i a_{k j}+ \eta_j a_{k i}=
  B_{,j} a_{k \alpha}a^\alpha_i + B_{,i} a_{k\alpha}a^\alpha_j
\end{gather}

\begin{rem}
  \label{mu}
  For further use let us note that if $B=\const$, i.e., if
  $B_{,i}\equiv 0$, then the right-hand side of the last equation
  vanishes implying $\eta_i\equiv 0$. Then,
  \begin{equation}
    \mu_{,i}= 2B\lambda_i.
    \label{eq:mu1}
  \end{equation}
\end{rem}

Let us now alternate~(\ref{eq:xe}) with respect to $(i,k)$, rename
$j\longleftrightarrow k$ and add the result to~(\ref{eq:xe}). After this manipulation
only the terms that  are symmetric with respect to  $(j,k)$ remain, and we obtain

\begin{gather}
  \label{eq:xejk}
  \xi_i g_{jk}+\eta_i a_{jk}=B_{,i} a_{k \alpha}a^\alpha_j
\end{gather}

If $B_{,i}\neq 0$, equation~(\ref{eq:xejk}) implies that for certain
 functions $C$ and $D$
\begin{gather}
  \label{eq:CD}
  C g_{jk}+D a_{jk}=a_{k \alpha}a^\alpha_j
\end{gather}

Let us now calculate $\nabla_k(a_{i \alpha}a^{\alpha}_j)$:
\begin{multline}
  \nabla_k(a_{i\alpha}a^{\alpha}_j)=
  a_{i\alpha,k}a^\alpha_j+a_{j\alpha,k}a^\alpha_i=\\
  = \Jij(\lambda_\i a_{\j k} + \lambda_\j a_{\i k}+\lambda_\alpha
  a^\alpha_\i g_{\j k} + \lambda_\alpha a^\alpha_\j g_{\i
    k})\stackrel{(\ref{eq:CD})}{=}\\
  =C_{,k} g_{ij}+D_{,k} a_{ij} + D \Jij(\lambda_\i g_{\j k} + \lambda_\j g_{\i k})
\end{multline}

Setting $s_i:=\lambda_\alpha a^\alpha_\i -D\lambda_i$, we obtain

\begin{gather}
  \label{eq:cd2}
  \Jij(\lambda_\i a_{\j k} + \lambda_\j a_{\i k}+s_i g_{\j k} + s_j
  g_{\i k}) -C_{,k} g_{\i\j}-D_{,k} a_{\i\j}=0
\end{gather}

To simplify this equation consider the action of the operator $\Jkj$
on it. After applying the properties of the complex structure, the equation
takes the form
\begin{gather}
  \Jij\left(\lambda_\i a_{\j k} + s_\i g_{\j k}\right)-
  \Jkj\left(\frac{C_{,\k}}{2}g_{i \j}+\frac{D_{,\k}}{2}a_{i \j}\right)=0.
\end{gather}
Alternating with respect to $(i,k)$ and collecting the terms yields
\begin{gather}
  \label{eq:Jijlag}
  \Jij\left[a_{\j k}\left(\lambda_\i + \frac{D_{,\i}}{2}\right) + g_{\j k}\left(s_\i +
    \frac{C_{,\i}}{2}\right)\right]
  -
  \Jkj \left[a_{i \j}\left(\lambda_\k + \frac{D_{,\k}}{2}\right) + g_{i \j}\left(s_\i +
    \frac{C_{,\i}}{2}\right)\right]
  =0
\end{gather}

After denoting
\begin{gather}
  \tau_i=s_i+\frac{C_{,i}}{2},\quad \bar\tau_i=J_{\ \ i}^\i\tau_\i,\quad
  g_{j \bar{k}}=J_{\ \ k}^\k g_{j \k}\\
  \nu_i=\lambda_i+\frac{D_{,i}}{2},\quad \bar\nu_i=J_{\ \ i}^\i\nu_\i,\quad
  a_{j \bar{k}}=J_{\ \ k}^\k a_{j \k},
\end{gather}
 equation~(\ref{eq:Jijlag}) reads
\begin{gather}
  \label{eq:tau_nu}
  (\tau_i g_{jk}-\tau_k g_{i j})-(\bar\tau_i g_{j \bar k} -
  \bar\tau_k g_{j \bar i})
  +(\nu_i a_{jk}-\nu_k a_{i j})-(\bar\nu_i a_{j \bar k} -
  \bar\nu_k a_{j \bar i})=0
\end{gather}
Let us now contract this equation with a certain vector field
$\xi^j$. We obtain
\begin{gather}
  \label{eq:tauxi}
  (\tau_i \xi_k-\tau_k \xi_i)-(\bar\tau_i \bar\xi_k -
  \bar\tau_k \bar\xi_i)=(\nu_i \eta_k-\nu_k \eta_i)-(\bar\nu_i \bar\eta_k -
  \bar\nu_k \bar\eta_i)
\end{gather}
where $\bar\xi_i=J_{\ \ i}^\i\xi_\i$, $\eta_i=-a_{ij}\xi^j$ and
$\bar\eta_i=J_{\ \ i}^\i\eta_\i$.

\par If the vectors $\tau_i$, $\xi_i$, $\bar\tau_i$ and $\bar\xi_i$ are
linearly independent, this equation implies that the $4$-dimensional space
$l(\tau,\xi)$ spanned over $\{\tau_i,\xi_i,\bar\tau_i,\bar\xi_i\}$
coincides with $l(\nu, \eta)$ spanned over $\{\nu_i,\eta_i,\bar\nu_i,
\bar\eta_i\}$. Indeed, these spaces are determined as the orthogonal
complements to the kernels of the corresponding 2-forms

\begin{gather*}
  Ker(\tau,\xi)=\{u^i\mid \left((\tau_i \xi_k-\tau_k \xi_i)-(\bar\tau_i
  \bar\xi_k - \bar\tau_k \bar\xi_i)\right)u^i x^k=0 \mbox{ for every
    $x^k$}\}\\
  Ker(\nu,\eta)=\{u^i\mid \left((\nu_i \eta_k-\nu_k \eta_i)-(\bar\nu_i
  \bar\eta_k -\bar\nu_k \bar\eta_i)\right)u^i x^k=0 \mbox{ for every
    $x^k$}\}
\end{gather*}
Since by~(\ref{eq:tauxi}) the forms are equal, the subspaces are equal as well.

\par If $\dim M\geq 6$,  there exist two vectors $\stackrel{1}{\xi^j}$
and $\stackrel{2}{\xi^j}$ such that
$\{\tau_i,\stackrel{1}{\xi_i},\stackrel{2}{\xi_i},
\bar\tau_i,\stackrel{1}{\bar\xi_i},\stackrel{2}{\bar\xi_i}\}$ are
linearly independent. Then $l(\tau,\stackrel{1}{\xi})$ and
$l(\tau,\stackrel{2}{\xi})$ intersect along   the $2$-dimensional subspace
spanned by the vectors $\{\tau_i,\bar\tau_i\}$. The corresponding vectors
$\stackrel{1}{\eta}$ and $\stackrel{2}{\eta}$ determine spaces
$l(\nu,\stackrel{1}{\eta})$ and $l(\nu,\stackrel{2}{\eta})$ which
intersect along  the subspace spanned by the vectors $\{\nu_i,\bar\nu_i\})$. Since
the $4$-dimensional spaces are pairwise equal, one obtains

\begin{gather}
  span\{\tau_i,\bar\tau_i\}=span\{\nu_i,\bar\nu_i\}
\end{gather}

\par Then, for certain functions $p,q$ we have
$$\tau_i=p\nu_i +q\bar\nu_i,\quad \bar\tau_i=p\bar\nu_i -q\nu_i.$$
Let us now substitute this in~(\ref{eq:tau_nu}). After collecting terms, we obtain
\begin{multline}
\label{eq:withoutomega} 
 \nu_i(p\,g_{jk}+q\,J^\k_{\ \ k} g_{j\k}+a_{jk})-
  \nu_k(p\,g_{ij}+q\,J^\i_{\ \ i} g_{\i j}+a_{i j})
  =\\
  \bar\nu_i(p\,J^\k_{\ \ k} g_{j\k}-q\,g_{jk}+J^\k_{\ \ k} a_{jk})-
  \bar\nu_k(p\,J^\i_{\ \ i}g_{\i j}-q\,g_{i j}+J^\i_{\ \ i}a_{\i j}).
\end{multline}

Defining
\begin{gather}
  \label{eq:omega}
  \omega_{j k}=p\,g_{jk}+q\,J^\k_{\ \ k} g_{j\k}+a_{jk},\\
  \omega_{j \bar{k}}=p\,J^\k_{\ \ k} g_{j\k}-q\,g_{jk}+J^\k_{\ \ k} a_{jk},
\end{gather}
 we can rewrite equation (\ref{eq:withoutomega}) in the form

\begin{gather}
  \nu_i \omega_{j k}-\nu_k \omega_{j i}
  =\bar\nu_i \omega_{j \bar k}- \bar\nu_k \omega_{j \bar{\imath}}
\end{gather}

This equation has the same structure as~(\ref{eq:tau_nu}), but with
a non-symmetric, hermitian bilinear form $\omega_{j k}$. One can
easily see that it holds if and only if

\begin{gather}
  \omega_{j k}=\alpha_j \nu_k + J^\j_{\ \ j} J_{\ \ k}^\k \alpha_\j \nu_\k
\end{gather}

for some covector $\alpha_j$.

\par Substituting $\omega$ in~(\ref{eq:omega}) and alternating the result, we obtain 

$$2 q J^\k_{\ \ k} g_{j \k}=\alpha_j \nu_k - \alpha_k \nu_j +J_{\ \ j}^\j J_{\ \ k}^\k
(\alpha_\j \nu_\k - \alpha_\k \nu_\j).$$

Let us now consider this equation as an equality between two bilinear
forms. The rank of the right-hand side is not greater then 4, while the
left hand side is nondegenerate unless $q\neq 0$. Since $\dim M\geq 6$
we have $q=0$ and $\omega_{j k}$ is symmetric
by~(\ref{eq:omega}). Thus,

\begin{gather}
  \omega_{j k}=\alpha (\nu_j \nu_k + J_{\ \ j}^\j J^\k_{\ \ k}\nu_\j \nu_\k),
\end{gather}

where $\alpha$ is a scalar function. It immediately follows that
(after renaming of variables)

\begin{gather}
  \label{eq:avbg}
  a_{ij}=p(u_i u_j+J^\i_{\ \ i} J^\j_{\ \ j} u_\i u_\j) + q g_{i j},
\end{gather}
where $p$, $q$ --- are certain functions and $u_i$ is a covariant
vector field.

\par We have shown that if in a neighborhood of some point there are
two linearly independent solutions of the extended system with
non-constant $B$, then each solution has the special
form~(\ref{eq:avbg}).

\par Now we would like to show that the function $q$, corresponding to
a solution $a_{ij}$ as was given in equation~(\ref{eq:avbg}), is a
constant. In order to do  this, take an arbitrary
$U^i\in\mbox{span}\{u^i,\bar{u}^i\}^{\perp}$. Contracting~\eqref{eq:avbg}
with $U^i$ we see that

\begin{gather}
  a_{i\alpha}U^{\alpha}= q U_{i}\nonumber
\end{gather}

Hence all vectors, orthogonal to $u$ and $\bar{u}$, correspond to the
eigenvalue $q$ of $a_{j}^{i}=g^{i\alpha}a_{\alpha j}$. Taking the derivative of the equation
above and inserting equation (\ref{eq:hpr}) yields

\begin{gather}
  \lambda_{i}U_{k}+\lambda_{\alpha}U^{\alpha}g_{ik}-\bar{\lambda}_{i}\bar{U}_{k}-\bar{\lambda}_{\alpha}U^{\alpha}J_{ik}+a_{i\alpha}U^{\alpha}_{,k}= q_{,k} U_{i}+q U_{i,k}\nonumber
\end{gather}

Contracting this equation with $U^{i}$ gives
\begin{gather}
  2\lambda_{\alpha}U^{\alpha}U_{k}-2\bar{\lambda}_{\alpha}U^{\alpha}\bar{U}_{k}=
  q_{,k} U_{\alpha}U^{\alpha}.
  \label{eq:lincomb3}
\end{gather}

Thus, $q_{,k}\in\mbox{span}\{U_{k},\bar{U}_{k}\}$ unless $U_\alpha U^\alpha=0$.

\par Given any vector $U^{i}\in\mbox{span}\{u,\bar{u}\}^{\perp}$, such
that $U_{\alpha}U^{\alpha}\neq 0$, we can construct a second vector
$W^{i}\in\mbox{span}\{u,\bar{u}\}^{\perp}$ such that
$W_{\alpha}W^{\alpha}\neq 0$ and
${\mbox{span}\{U^{i},\bar{U}^{i}\}}\cap{\mbox{span}\{W^{i},\bar{W}^{i}\}=\{0\}}$.
In this case, using equation~(\ref{eq:lincomb3}) for $U^{i}$ and
$W^{i}$, we obtain that $q$ is a constant (because $q_k\in\mbox{span}\{U_{k},\bar{U}_{k}\}\cap\mbox{span}\{W_{k},\bar{W}_{k}\}=\{\vec{0}\}$). It remains to show
that such a vector $U^{i}$ exists. Assuming each vector
$U^{i}\in\mbox{span}\{u,\bar{u}\}^{\perp}$ satisfies
$U_{\alpha}U^{\alpha}=0$, we obtain that $U_{\alpha}W^{\alpha}=0$ for
all $U^{i},W^{i}\in\mbox{span}\{u,\bar{u}\}^{\perp}$. Since $\mbox{dim}\,M\geq 6$, this means that $\mbox{dim}\,\mbox{span}\{u,\bar{u}\}=\mbox{dim}\,((\mbox{span}\{u,\bar{u}\})^{\perp})^{\perp}\geq 4$ which is a contradiction.

 Using that $q$ is a constant, we can substract the trivial solution $qg_{ij}$ from $a_{ij}$ and
include the function $p$ in the vector field $u_i$. In other words, without loss of generality, 
$a_{ij}$ is  given by

\begin{gather}
  a_{ij}=u_i u_j+J^\i_{\ \ i} J^\j_{\ \ j} u_\i u_\j.\nonumber
\end{gather}

Note that  $u^{j}$ is an eigenvector of $a^{i}_{j}$ as well. If
the corresponding  eigenvalue is a constant,  all eigenvalues of   $a^{i}_j$ are constant.  Hence,  the trace of $a^{i}_j$ is constant, and the $1$-form
$\lambda_{i}=\tfrac{1}{4} (a^k_k)_{,i}$ is identically zero. Inserting
$\lambda_{i}\equiv 0$ in equation (\ref{eq:mg+Ba1}), we see that

\begin{align}
  0=\mu g_{ij}+Ba_{ij}\nonumber
\end{align}

By Lemma~\ref{lem:globind}, $\mu=B=0$, since $g_{ij}$ and $a_{ij}$ are assumed to be linearly independent. We see that in this case $B=\const$ as we claim. 

\par Now consider the case when the eigenvalue  corresponding to the
eigenvector $u^{i}$ is not  constant. We obtain that
$\mbox{span}\{\lambda_{i},\bar{\lambda}_{i}\}=\mbox{span}\{u_{i},\bar{u}_{i}\}$,
since $\lambda_{i}$ and $\bar{\lambda}_{i}$ are contained in the sum
of the eigenspaces, corresponding to the non-constant
eigenvalues. Consider the second solution

\begin{gather}
  A_{ij}=v_i v_j+J^\i_{\ \ i} J^\j_{\ \ j} v_\i v_\j\nonumber
\end{gather}

of the extended system, such that
% the corresponding function $B$ is not a constant and
$a_{ij}$, $A_{ij}$, $g_{ij}$ are linearly independent. By
$\Lambda_{i}$, we denote the $1$-form corresponding to $A_{ij}$. The sum

\begin{gather}
  a_{ij}+A_{ij}=u_i u_j+J^\i_{\ \ i} J^\j_{\ \ j} u_\i u_\j+v_i v_j+J^\i_{\ \ i} J^\j_{\ \ j} v_\i v_\j \nonumber
\end{gather}
is again a solution of equation (\ref{eq:hpr}) and hence, can be written as
\begin{gather}
  a_{ij}+A_{ij}=w_i w_j+J^\i_{\ \ i} J^\j_{\ \ j} w_\i w_\j+Q g_{ij}\nonumber
\end{gather}

Comparing the last two equations, we see that
\begin{gather}
  Q g_{ij}=u_i u_j+\bar{u}_i \bar{u}_j+v_i v_j+\bar{v}_i \bar{v}_j-w_i w_j-\bar{w}_i \bar{w}_j.\nonumber
\end{gather}
 Since
$\mbox{span}\{\lambda_{i},\bar{\lambda}_{i}\}=\mbox{span}\{u_{i},\bar{u}_{i}\},
\mbox{span}\{\Lambda_{i},\bar{\Lambda}_{i}\}=\mbox{span}\{v_{i},\bar{v}_{i}\}$
and
$\mbox{span}\{w_{i},\bar{w}_{i}\}=\mbox{span}\{\lambda_{i}+\Lambda_{i},\bar{\lambda}_{i}+\bar{\Lambda}_{i}\}$,
the right-hand side has rank at most $4$ and therefore, $Q\equiv 0$. Let us rewrite the last equation in the form
\begin{gather}
  w_i w_j+\bar{w}_i \bar{w}_j=u_i u_j+\bar{u}_i \bar{u}_j+v_i v_j+\bar{v}_i \bar{v}_j\nonumber
\end{gather}

Since the left hand side has rank 2, $u_{i}$, $\bar{u}_{i}$, $v_{i}$ and
$\bar{v}_{i}$ are linearly dependent and the intersection
$\mbox{span}\{u_i,\bar{u}_i\}\cap\mbox{span}\{v_i,\bar{v}_i\}$ is
non-empty. Since it is also $J$-invariant, we obtain that
$\mbox{span}\{u_i,\bar{u}_i\}=\mbox{span}\{v_i,\bar{v}_i\}$. Thus,  $v_{i}=\alpha u_i+\beta \bar{u}_i$, for some real constants
$\alpha,\beta$. It follows, that $\bar{v}_{i}=\alpha \bar{u}_i-\beta
u_i$ and we obtain
\begin{align}
  v_i v_j+\bar{v}_i \bar{v}_j&=(\alpha u_i+\beta \bar{u}_i)(\alpha
  u_j+\beta \bar{u}_j)+(\alpha \bar{u}_i-\beta u_i)(\alpha
  \bar{u}_j-\beta u_j)\nonumber\\
  &=(\alpha^{2}+\beta^{2})(u_iu_j+\bar{u}_i\bar{u}_j)\nonumber
\end{align}

Inserting this in the original formulas for $a_{ij}$ and $A_{ij}$, we
see that $a_{ij}=\const\cdot A_{ij}$. We obtain  a contradiction
to the  assumption that $a_{ij}$ and $A_{ij}$ are linearly
independent.  Lemma~\ref{lem:ext2} is proved under the assumption  $\dim
M\geq 6$.

\subsubsection{Proof of Lemma~\ref{lem:ext2} in case $\dim M=4$.}
\label{Bconstdim4}

\begin{lem}
  \label{lem:dimM4}
  Let $(M^{2n},g,J)$ be a Kähler manifold of dimension $2n=4$ and assume that the degree of mobility of the metric $g$ is $\geq3$. Then $g$ has constant holomorphic sectional curvature $-4B$,
  where $B$ is defined by  equation (\ref{eq:mg+Ba1}). In particular, $B$ is a constant.
\end{lem}

\begin{rem}
  As we see, Lemma~\ref{lem:dimM4} contains an extra statement: not
  only $B=\const$, but also the metric $g$ has constant holomorphic
  sectional curvature. This result was actually unexpected. Indeed,
  the analog of dimension $4$ in the theory of projectively equivalent
  metrics is $2$, and in dimension $2$ there exist metrics of
  non-constant sectional curvature admitting $4$-parametric family of
  projectively equivalent metrics.
\end{rem}

\begin{proof}
  We will work in a small neighborhood of the point $p\in M$, such that
  there exist three solutions $g_{i j}$, $a_{i j}$ and $A_{i j}$ of
  equation (\ref{eq:hpr}), linearly independent at $p$.

  % ystem~(\ref{eq:hpr}) every its solution $a_{ij}$ satisfies the
  % following equation:
  % \begin{align}
  %   a_{i\alpha}R^{\alpha}_{jkl}+a_{j\alpha}R^{\alpha}_{ikl}=\lambda_{i,l}g_{jk}
  %   -\lambda_{i,k}g_{jl}+\lambda_{j,l}g_{ik}-\lambda_{j,k}g_{il}-\bar{\lambda}_{i,l}J_{jk}
  %   +\bar{\lambda}_{i,k}J_{jl}-\bar{\lambda}_{j,l}J_{ik}+\bar{\lambda}_{j,k}J_{il}\nonumber
  % \end{align}

  Using equation (\ref{eq:mg+Ba1}), we substitute $\lambda_{i,j}$ in equation~\eqref{integr} to obtain
\begin{equation}  
a_{i\alpha}R^{\alpha}_{jkl}+a_{j\alpha}R^{\alpha}_{ikl}=
  -4B (a_{i\alpha}K^{\alpha}_{jkl}+a_{j\alpha}K^{\alpha}_{ikl}),\label{eq:com1}
\end{equation}
  where $K$ is the algebraic curvature tensor of constant holomorphic sectional curvature equal to $1$, namely

  \begin{align}
    K^{\alpha}_{jkl}=\frac{1}{4}(\delta^{\alpha}_{k}g_{jl}-\delta^{\alpha}_{l}g_{jk}+
    J^{\alpha}_{\ \ k}J_{jl}-J^{\alpha}_{\ \ l}J_{jk}+2J^{\alpha}_{\ \ j}J_{kl}).\nonumber
  \end{align}

  Let us define the $(1,3)$-tensor $G^i_{jkl}=R^i_{jkl}+4 B K^i_{jkl}$. This new 
  tensor has the same algebraic symmetries as the Riemannian
  curvature tensor $R$ (including the Bianci identity), in particular, it commutes with the complex
  structure $J$:
  \begin{gather}
    \label{eq:Gsym}
    G_{ij kl}=-G_{ji kl}\ ,\\
    G_{ij kl}=G_{kl ij}\ ,
    G^i_{\alpha kl}J_{\ \ j}^{\alpha}=J^i_{\ \ \alpha}G^\alpha_{jkl}
  \end{gather}
  In addition, from equation (\ref{eq:com1}) it follows, that $G^{i}_{jkl}$ satisfies 
  \begin{gather}
    \label{eq:Ga}
    a_{i\alpha}G^\alpha_{jkl}+a_{\alpha j} G^\alpha_{ikl}=0
  \end{gather}
for each solution $(a_{ij},\lambda_{i})$ of equation (\ref{eq:hpr}), $a_{ij}\neq\mbox{const}\cdot g_{ij}$.\\  
Our goal is to show that $G^{i}_{jkl}\equiv 0$.\\
  For an arbitrary skew-symmetric $(2,0)$-tensor $\omega^{kl}$
  consider the linear operator $$G(\omega)^i_j\mathrel{\mathop:}=G^i_{jkl}\omega^{kl}.$$

  \par Since $g$ is hermitian, there exists a basis in $T_p
  M$ such that the matrices of $g$ and $J$ are given by
  $$g=
  \begin{pmatrix}
    1 &&& \\
    &1&& \\
    &&\varepsilon& \\
    &&&\varepsilon
  \end{pmatrix}, \qquad
  J=\begin{pmatrix}
    &-1&& \\
    1&&& \\
    &&&-1 \\
    &&1&
  \end{pmatrix}$$
  where $\varepsilon=\pm 1$ depending on the signature of $g$. Fixing this basis, we will
  work in matrix notation. Since $g$ is non-trivial, it is important
  to note that letters $J$, $a$ and $G(\omega)$ correspond to matrices
  of linear operators, i.e. $(1,1)$-tensors. By $g$ we denote the matrix
  of the $(0,2)$-form $g_{i j}$.
  \par All matrices we are working with commute with the complex structure
  $J$. It is a well-known fact (that can be checked by direct
  calculation) that matrices, commuting with the complex structure, are
  ``complex'' in the sence that they have the form
  \begin{gather}
    \begin{pmatrix}
      \alpha_1 & \beta_1&\alpha_2 & \beta_2\\
      -\beta_1 & \alpha_1&-\beta_2 & \alpha_2 \\
      \alpha_3 & \beta_3&\alpha_4 & \beta_4\\
      -\beta_3 & \alpha_3&-\beta_4 & \alpha_4
    \end{pmatrix}
  \end{gather}

  Using this form one can define the nondegenerate $\R$-linear mapping $\psi$
  \begin{gather}
    \psi:\{Q\in \textrm{Mat}(4,4,\mathbb R)\mid QJ=JQ\}\to \textrm{Mat}(2,2,\mathbb C)
  \end{gather}
  given by the formula
  $$\psi\left(
    \begin{pmatrix}
      \alpha_1 & \beta_1&\alpha_2 & \beta_2\\
      -\beta_1 & \alpha_1&-\beta_2 & \alpha_2 \\
      \alpha_3 & \beta_3&\alpha_4 & \beta_4\\
      -\beta_3 & \alpha_3&-\beta_4 & \alpha_4
    \end{pmatrix}
  \right)=
    \begin{pmatrix}
      \alpha_1 + i \beta_1&\alpha_2 + i \beta_2\\
      \alpha_3 + i \beta_3&\alpha_4 + i \beta_4\\
    \end{pmatrix}
  .$$

  It is easy to check that $\psi(Q_1 Q_2)=\psi(Q_1)\psi(Q_2)$ and
  $\psi(Q^T)=\overline{\psi(Q)}^T, $ where ``$\overline{\hspace{2ex}}$'' denotes the complex conjugation.   Moreover,  $\psi(J)=i\cdot \Id$ and
  $\psi(g)=
  \begin{pmatrix}
    1&0\\
    0&\varepsilon\\
  \end{pmatrix}$.

  To simplify the notation we will identify a matrix with its image
  under the mapping $\psi$, for example $a$ and $\psi(a)$ are identified, as well as $g$ and
  $\psi(g)$.

  \par Since $a_{ij}$ is symmetric, it satisfies the equation
  \begin{gather}
    g a=\overline{(g a)}^T
  \end{gather}

  Thus, there exist real numbers $\alpha$, $\beta$ and a complex number $Z$ such that
  \begin{gather}
    a=
    \begin{pmatrix}
      \alpha&Z\\
      \bar Z&\beta\\
    \end{pmatrix}.
  \end{gather}

  By assumptions there exist three solutions $a_{ij}$, $A_{ij}$,
  $g_{ij}$ which are linearly independent at the point. Then there
  exists a nontrivial (i.e., $\ne c \cdot g$ at the point we are
  working in) solution such that $\alpha= \beta =0$. Without loss of
  generality we think that the solution $a_{ij}$ has $\alpha= \beta
  =0$ and $Z\neq 0$, i.e.

  \begin{gather}
    \label{eq:aZ}
    a=
    \begin{pmatrix}
      0&Z\\
      \bar Z&0\\
    \end{pmatrix}
  \end{gather}

  Consider now the restrictions that  equations~\eqref{eq:Gsym}
  and~\eqref{eq:Ga} impose on the complex form of $G(\omega)$. Since
  $G(\omega)_{ij}$ is skew-symmetric
  $$g G(\omega)=-\overline{(g G(\omega))}^T.$$
  Thus, $G(\omega)$ has the form
  \begin{gather}
    \label{eq:Gomega}
    G(\omega)=
    \begin{pmatrix}
      i\alpha&W\\
      -\overline W&i\beta\\
    \end{pmatrix}
  \end{gather}
  for certain real numbers $\alpha$, $\beta$ and a complex number $W$.
  The last condition we have to make use of is
  $$a G(\omega) = G(\omega) a$$

  \par Since $a$ is simple (moreover, has different eigenvalues) every
  matrix that commutes with $a$ is a polynomial of $a$. (Recall that
  the matrix $a$ in our convention corresponds to the $(1,1)$-tensor
  $a^i_j$.)  Thus, $G(\omega)=C\cdot a + D \cdot \Id$ for certain
  complex numbers $C$ and $D$. Using the explicit form of $a$ and
  $G(\omega)$ (see~\eqref{eq:aZ} and~\eqref{eq:Gomega}) we obtain
  \begin{gather}
    \begin{pmatrix}
      i\alpha&W\\
      -\overline W&i\beta\\
    \end{pmatrix}
    = C
    \begin{pmatrix}
      0&Z\\
      \bar Z&0\\
    \end{pmatrix}
    + D
    \begin{pmatrix}
      1&0\\
      0&1\\
    \end{pmatrix}
  \end{gather}
  which implies that both $D=i \alpha=i \beta$ and
  $C=\frac{W}{Z}=-\overline{\left(\frac{W}{Z}\right)}$ are purely
  imaginary.  Finally we obtain
  \begin{gather}
    \label{eq:Gaid}
    G(\omega)=i\cdot  c \cdot  a + i\cdot d \cdot \Id
  \end{gather}
  with real coefficients $c$, $d$. If we assume $c\neq 0$, then $G(\omega)$
  has different eigenvalues. Thus, $G(\omega)$ is simple. Let us consider
  another solution $A$ of equation~\eqref{eq:hpr}. Since it commutes
  with the simple matrix $G(\omega)$ it is a polynomial of $G(\omega)$:
  \begin{gather}
    \label{eq:AG}
    A= \tau G(\omega) + \nu \,\Id
  \end{gather}
  Subsituting the explicit form of $A=
  \begin{pmatrix}
    \alpha_A&Z_A\\
    \bar Z_A&\beta_A\\
  \end{pmatrix}
  $ we obtain
  \begin{gather}
    \begin{pmatrix}
      \alpha_A&Z_A\\
      \overline Z_A&\beta_A\\
    \end{pmatrix}=
    \tau
    \begin{pmatrix}
      i\alpha&W\\
      -\overline W&i\beta\\
    \end{pmatrix}
    + \nu
    \begin{pmatrix}
      1&0\\
      0&1\\
    \end{pmatrix}
  \end{gather}
  which implies that $\tau=i\, t$ is purely imaginary and $\nu$ is
  real. Therefore, equation~\eqref{eq:AG} implies that all solutions
  of~\eqref{eq:mg+Ba1} are contained in the $2$-dimensional space $i t
  G(\omega) + \nu \,\Id$, which gives us the contradiction. Then, $c =
  0$. Thus, from~\eqref{eq:Gaid} we obtain that for every $\omega$
  the operator $G(\omega)$ is proportional to the complex structure: in the initial ``real'' notation, we obtain 

  \begin{gather}
    G(\omega)^i_j= d(\omega)J^i_{\ \ j}
  \end{gather}

  Since the left hand side is linear in $\omega^{kl}$,
  it follows that $d(\omega)=d_{kl}\omega^{kl}$ and hence
  $G^{i}_{jkl}\omega^{kl}=d_{kl}\omega^{kl}J^{i}_{\ \ j}$ implying
  $G_{ijkl}=d_{kl}J_{ij}$.  Using the symmetry
  relations~\eqref{eq:Gsym} for $G_{ijkl}$ we obtain
  $d_{kl}J_{ij}=d_{ij}J_{kl}$ and therefore $d_{kl}=c J_{kl}$ for some
  constant $c\neq 0$. Let us show that $G_{ijkl}=c J_{ij} J_{kl}$ does
  not satisfy the Bianci identity unless $c=0$.
  % Let us now consider the Bianchi identity for
  % $G_{ijkl}=c J_{ij} J_{kl}$.
  By direct computation we obtain
  $$0=G_{1234}+G_{1423}+G_{1342}=c(J_{12}J_{34}+J_{14}J_{32}+J_{13}J_{42})
  =c(1\cdot1+0\cdot0+0\cdot0)=c.$$
  Thus, $G^{i}_{jkl}\equiv 0$.

  \par Finally,
  % \par We have shown that at the neighborhood of every point such as
  % there exist three linearly independent solutions
  % of~\eqref{eq:mg+Ba1}
  $$0=G^{i}_{jkl}=R^{\alpha}_{jkl}+4BK^{\alpha}_{jkl},$$
  i.e. our metric $g$ has pointwise constant holomorphic curvature
  $-4B$ (at almost every point, and therefore at every point of $M$).
  Thus, $M$ has constant holomorphic sectional curvature (see for
  example~\cite[chapter $8$]{Kob}). Then $B$ is a constant and Lemma~\ref{lem:ext2} has been proved for $\mbox{dim}\,M=4$. 

\end{proof}

% Let us now recall that equation~\eqref{eq:Gaid}} holds for an
% arbitrary skew-symmetric $(2,0)$-tensor $\omega^{kl}$. Suppose that
% there exists $\omega'$ such as
% \begin{gather}
% G(\omega')=i\, c' \, a + i\, d' \, Id
% \end{gather}
% with real numbers

%\subsection{The constant $B$ and the function $\mu$ are globally defined}
\subsection{Last step in the proof of Theorem~\ref{thm:ext}}
\label{subsec:globalB}

Above,  we proved the following 

{\bf Statement.} {\it Let $(M^{2n},g,J)$ be a connected Kähler manifold of
  dimension $2n\geq 4$. Assume the degree of mobility $D$ of  $g$ is
  $\geq 3$. Then, for every solution $(a_{ij}, \lambda_i)$  of  \eqref{eq:hpr} such that   $a_{ij}\ne \const \cdot g_{ij}$, 
  almost
 every point of $M$ has a neighborhood such that in this neighborhood
  there exists an unique constant $B$ and a scalar function $\mu$ such
  that the ``extended'' system~(\ref{eq:mg+Ba}) holds.}

Indeed, the first equation of (\ref{eq:mg+Ba}) is  equation~\eqref{eq:hpr} and is
fulfilled everywhere. The second equation is fulfilled almost
everywhere by the results of the previous sections. Now, as we noted
in Remark~\ref{mu}, at every open set such that the second equation is
fulfilled, the third equation is fulfilled as well.

\begin{rem}
  The above {\bf statement } is visually close to
  Theorem~\ref{thm:ext}, the only difference is that in
  Theorem~\ref{thm:ext} the constant $B$ and the function $\mu$ are
  universal (i.e., do not depend on the neighborhood).  We will prove
  it in this section.
\end{rem}

First let us prove

\begin{lem}
  \label{lem:tanno}
  Assume that in every point of an open subset $U\subseteq M$ the
  extended system~(\ref{eq:mg+Ba}) holds (for a certain constant
  $B$). Then, in this neighborhood, the function $\lambda:=
  \tfrac{1}{4}a^i_i$ satisfies Tanno's equation
  \begin{gather}
    \label{tanno}
    \lambda_{,ijk}=B(2\lambda_{,k} g_{ij}+ \lambda_{,i} g_{j k} +
    \lambda_{,j} g_{i k} -
    \bar{\lambda}_{,i}J_{jk}-\bar{\lambda}_{,j}J_{ik}).
  \end{gather}
\end{lem}

\begin{rem}
Recall that the differential  of the function  $\lambda$  is precisely the covector $\lambda_i$ from \eqref{eq:a}, i.e., $\lambda_{,i}=\lambda_i, $ see the discussion after Theorem \ref{thm:mikes}.
\end{rem}
\begin{proof}
  If $B$ is a constant, the function $\mu$ is smooth as the
  coefficient of the proportionality of the nonvanishing smooth tensor
  $g_{ij}$ and the smooth tensor $(\lambda_{i,j} - Ba_{ij})$. 

   We take the covariant derivative of the second equation of the
  ``extended'' system and substitute the first  and the third equations  
  inside. In view of $\lambda_i= \lambda_{,i}$, we obtain 
  \begin{align*}
    \lambda_{,ijk}&=\mu_{,k} \cdot g_{ij}+Ba_{ij,k}=
    2B\lambda_k \cdot g_{ij} + B(\lambda_i g_{j k} + \lambda_j g_{i k}
    - \bar{\lambda}_{i}J_{jk}-\bar{\lambda}_{j}J_{ik})\\
    &=
    B(2\lambda_k \cdot g_{ij}+ \lambda_i g_{j k} + \lambda_j g_{i k}
    - \bar{\lambda}_{i}J_{jk}-\bar{\lambda}_{j}J_{ik}).
  \end{align*}
\end{proof}

\par Now let us prove that the constant $B$ is universal. It is
sufficient to prove this in a neighborhood $W(q)$ of an
\emph{arbitrary} point $q$. Indeed, every continuous curve $c:[t_0,
t_1]\to M^{2n}$ lies in finite number of such neighborhoods $W$. Since  
the constants $B$ for two such intersected neighborhoods  must coincide,  the value of $B$ at the point $c(t_0) $ equals the value
of $B$ at $c(t_1). $ Since the manifold is assumed to be connected,
the constant $B$ is therefore universal, i.e., is the same for all neighborhoods.

\par Let $W\subseteq M$ be a sufficiently small neighborhood. Without loss of generality we can 
 assume that $W$ is geodesically convex, that is,  every two points $p, \tilde p \in W$ can be
connected by a unique geodesic segment  lying in $ W$.

\par We want to show that each two open sets contained in $W$ such
that they are as in the {\bf statement} above have the same constant
$B$. Let $U,\tilde{U}\subseteq W$ be nonempty open sets such that in
these sets the extended equations (\ref{eq:mg+Ba}) are satisfied with constants $B$ for
$U$ and $\tilde{B}$ for $\tilde{U}$.

\par
\begin{figure}
  \centering
 \includegraphics{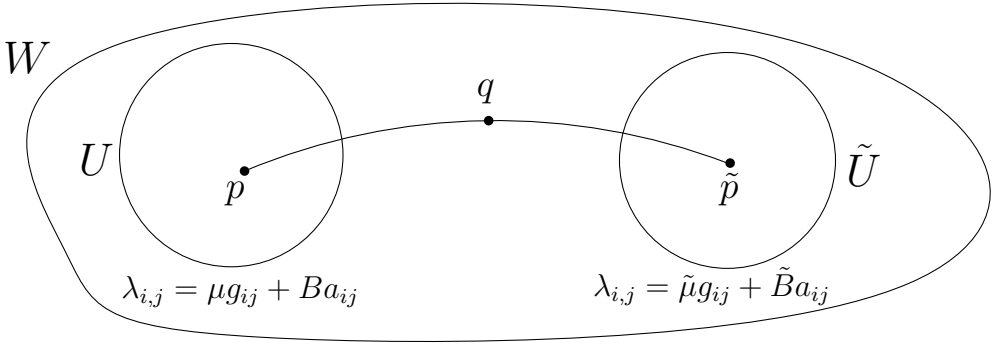} %\includegraphics{graph/p_q_tilde_p_new}
  \caption{There exists $q$ on $\gamma_{p,\tilde p}$ such that $\lambda_i=0$ at $q$.}
  \label{fig:p_q_tilde_p}
\end{figure}

 We assume $B\neq\tilde{B}$. We take a point $p\in U$ and connect
this point with every point $\tilde{p}\in\tilde{U}$ by a geodesic
$\gamma_{p,\tilde{p}}:[0,1]\rightarrow W$, $\gamma_{p,\tilde{p}}(0)=
p$, $\gamma_{p,\tilde{p}}(1)=\tilde p$ (see
Fig.~\ref{fig:p_q_tilde_p}).

Let us show that $\gamma_{p,\tilde{p}}$ contains a point $q$ such that
$\lambda_i= 0$ at $q$. Indeed, contracting equation~\eqref{tanno} with
$g^{ij}$ we obtain
\begin{align}
  \Delta\lambda_{,k}=4B(n+1)\lambda_{k}.
  \label{tanno1}
\end{align}

If $\lambda_i \ne 0$ at all points of the geodesic $\gamma_{p,\tilde
  p}$, we can find a vector field $\xi^i$ in some neighborhood
$U(\gamma_{p,\tilde p})$ of the geodesic $\gamma_{p,\tilde p}$ such
that $\lambda_i\xi^i\ne 0$ at all points of this neighborhood
$U(\gamma_{p,\tilde p})$.  Then, the function
\begin{equation}\label{81} 
  \tfrac{\Delta\lambda_{,k}\xi^{k}}{4(n+1)\lambda_{k}\xi^{k}}
\end{equation}
is well defined and smooth in $U(\gamma_{p,\tilde p})$. Comparing \eqref{tanno1} with \eqref{81}, we see that   in a
neighborhood of almost every point it is equal to the constant $B$ in
this neighborhood, so it is constant on $U(\gamma_{p,\tilde
  p})$. Then, $ B=\tilde B $ which contradicts our
assumption. Finally, there exists a point $q$ of the geodesic
$\gamma_{p, \tilde p}$ such that $\lambda_i=0$ at $q$.

\par By Corollary~\ref{cor:killing}, $\bar{\lambda}_{i}$ is a Killing
vector field. Then, the function $\dot \gamma_{p, \tilde p}^i \bar
\lambda_i$ is constant on the geodesic $\gamma_{p,\tilde p}$. Since it
vanishes at $q$, it vanishes at all other points of $\gamma_{p,\tilde
  p}$, in particular we have that at the point $p=\gamma_{p, \tilde
  p}(0)$ the vector $\bar\lambda^i$ is orthogonal to $\dot \gamma_{p,
  \tilde p}^i(0)$.

\begin{figure}[ht]
  \centering
  \includegraphics{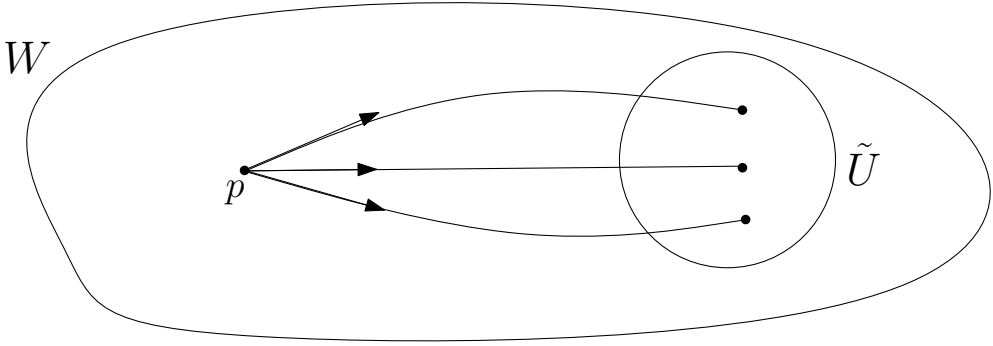}
  \caption{$\bar\lambda^i$ at $p$ is orthogonal to every $\dot
    \gamma^i_{p,\tilde p}$, implying $\lambda^i\equiv 0$.}
  \label{fig:lambda_at_p}
\end{figure}

The same is true for every geodesic connecting the point $p$
with any other point of $\tilde U$. Then, the vector $\bar \lambda^i $
at $p$ is orthogonal to many vectors (to all initial vectors of the
geodesics starting from $p$ and containing at least one point of
$\tilde U$); thus $\lambda_i=0$ at $p$ (see
Fig.~\ref{fig:lambda_at_p}).

\par Replacing the point $p$ by any other point of the neighborhood
$U$, we obtain that $\lambda_i=0$ at all points of $U$. By Corollary
\ref{cor:killing2}, $\lambda_i\equiv 0$ on the whole
manifold. Substituting $\lambda_i\equiv 0$ in the extended system, and
using that $g_{ij}$ is not proportional to $a_{ij}$, we see that
$B=0$ (at almost  all points of manifold).

\par Thus, the constant $B$ is universal on the whole connected
manifolds. Theorem~\ref{thm:ext} is proved.

\section{ The case $B=0$}
\label{sec:B=0}
By Corollary \ref{cor:Bind}, we already now that the global constant $B$, arising in the extended system (\ref{eq:mg+Ba}), does not depend on the solutions $(a_{ij},\lambda_{i})$ of equation (\ref{eq:hpr}). In this section we want to investigate the case when $B=0$. Our goal is to prove the following
\begin{thm}
  \label{thm:B=0}
  Let $(M^{2n}, g, J)$ be a closed connected Kähler manifold of dimension $2n\geq 4$ and of degree of mobility $\geq 3$. Suppose the constant $B$ in the system (\ref{eq:mg+Ba}) is zero, then $\lambda_{i}\equiv 0$ on the whole  $M$ for each solution $(a_{ij},\lambda_i)$ of equation (\ref{eq:hpr}).
  
In particular, every  metric $\bar{g}$, $h$-projectively equivalent to $g$, is already affinely equivalent to $g$.
\end{thm}

\begin{proof}
  If $B=0$, then $\mu=\const$ by the third equation
  from~(\ref{eq:mg+Ba}), and the second equations reads $\lambda_{i,j}
  = \const \cdot g_{ij}$. Then, the hessian $\lambda_{i,j} $ of the
  function $\lambda:= \tfrac{1}{4} a_i^i$ is covariantly constant.

  \par Since the manifold is closed the function $\lambda$ has a
  minimum and a maximum. At a minimum, the Hessian must be
  non-negatively definite, and at a maximum it must be nonpositively
  definite. Therefore the Hessian is null, and $\lambda_i$ is
  covariantly constant. But as it vanishes at the extremal points, it
  vanishes everywhere. Thus, $\lambda_i \equiv  0$ as we claim. By Remark \ref{aff}, every  metric $\bar{g}$, $h$-projectively equivalent to $g$, is already affine equivalent to $g$ as we claim. 
\end{proof}

\section{ If $B\ne 0$, the metric $-B\cdot g $ is positively definite}
\label{sec:proj}
Now let us treat the case when the constant $B$ in the system (\ref{eq:mg+Ba}) is different from zero.
Let $(M^{2n}, g, J)$ be a connected Kähler manifold of dimension $2n\geq
4$.  Let $(a_{ij}, \lambda_i, \mu)$ be a solution of ~(\ref{eq:mg+Ba}).
Since $B\neq 0$,  we can replace $g$ by the metric $-B\cdot g$ (having
the same Levi-Civita connection with $g$).

\par Then, for every solution $(a_{ij}, \lambda_i, \mu)$ of the
system~(\ref{eq:mg+Ba}), the triple $(-B \cdot a_{ij}, \lambda_i,
-\tfrac{1}{B}\mu)$ is the solution of ~(\ref{eq:mg+Ba}) corresponding
to the metric $g':= -B \cdot g$ with the constant $B=-1$.  Indeed, the
Levi-Civita connections of $g$ and $g'$ coincide, so substituting $(-B
\cdot a_{ij}, \lambda_i, -\tfrac{1}{B}\mu, -Bg, -1)$ instead of $(
a_{ij}, \lambda_i, \mu, g, B)$ in the extended system gives the system
which is equivalent to the initial extended system.

\par Note that the mapping $(a_{ij}, \lambda_i, \mu)\mapsto (-B \cdot
a_{ij}, \lambda_i, -\tfrac{1}{B}\mu)$ is linear and bijective, so the
degrees of mobility of $g$ and ${-}Bg$ are equal. Thus, if $B\ne 0$,
in the proof of Theorem~\ref{thm:main}, without loss of generality we
can assume that $B=-1$.

\par The goal of this section is to prove the following

\begin{thm}
  \label{pos}
  Let $(M^{2n}, g, J)$ be a closed connected Kähler manifold of dimension $2n\geq 4$.
  Suppose $(a_{ij},\lambda_i,\mu)$ satisfies
  \begin{gather}
    \label{eq:system}
    \begin{array}{c}
      a_{i j,k}=\Jij(\lambda_\i g_{\j k} + \lambda_\j g_{\i k})\\
      \lambda_{i,j}=\mu g_{i j}-a_{i j},\\
      \mu_{,i}=-2\lambda_i
    \end{array}
  \end{gather}
and $\lambda_{i}\neq 0$ at least at one point. Then, the metric $g$ is positively definite.
\end{thm}

\begin{rem}
  The assumption that the manifold is closed is important -- one can
  construct examples of complete pseudo-Riemannian Kähler
   metrics  admitting nontrivial solutions
  $(a_{ij},\lambda_i,\mu)$. Simplest  examples are 
   pseudo-Riemannian  Kähler  manifolds  
    of  constant   holomorphic curvature $4$. 
    Examples  of nonconstant holomorphic 
    curvature also exist and can be constructed similar to  \cite[Example 3.1]{Cortes2009}.
    \end{rem}

We need the following 
\begin{lem}
  \label{lem:frobenius}
  Let $(a_{ij},\lambda_i,\mu)$ be a solution of the
  system~(\ref{eq:system}) such that $a_{ij}=0$, $\lambda_i=0$,
  $\mu=0$ at some point $p$ of the connected Kähler manifold
  $(M^{2n},g,J)$.

  \par  Then $a_{ij}\equiv 0$, $\lambda_i\equiv 0$, $\mu\equiv 0$ at all
  points of $M$. In particular, the degree of mobility is always
  finite.
\end{lem}

\begin{proof}
  The system \eqref{eq:system} is in the Frobenius form, i.e., the
  derivatives of the unknowns $a_{ij}, \lambda_i, \mu$ are expressed
  as (linear) functions of the unknowns:
  $$
  \left(\begin{array}{c}a_{ij,k}\\\lambda_{i,j}\\\mu_{,i}\end{array}\right)
  =F\left(\begin{array}{c}a_{ij}\\\lambda_{i}\\\mu\end{array}\right)\nonumber,
  $$
  and all linear systems in the Frobenius form have the property that
  the vanishing of the solution at one point implies the vanishing at all
  points.
\end{proof}
The rest of this section is dedicated to the proof of Theorem~\ref{pos}. 
Our first goal will be to show, that it is possible to choose one solution of the system (\ref{eq:system}) (under the assumptions of Theorem~\ref{pos})
such that the corresponding operator $a^i_j=g^{i\alpha}a_{\alpha j}$ has a clear and simple structure of eigenspaces and eigenvectors. 

\subsection{Matrix of the extended system} \label{99}
In order to find the special solution of (\ref{eq:system}) mentioned above, we rewrite a solution $(a_{ij},\lambda_{i},\mu)$ as a  $(1,1)$-tensor on the $(2n+2)$-dimensional manifold $\widehat{M}=\R^2\times M$ with coordinates $(\underbrace{x_+,x_-}_{\R^2},\underbrace{x_1,\dots,x_{2n}}_{M})$. For every solution $(a^i_j,\lambda_i,\mu)$ of the system (\ref{eq:system}), let us consider  the $(2n+2)\times(2n+2)$-matrix
\begin{gather}
  \label{eq:operator}
  L(a,\lambda,\mu)=\left(
    \begin{array}{cc|ccc}
      \mu&0&\lambda_{1}&\dots&\lambda_{2n}\\
      0&\mu&\bar{\lambda}_{1}&\dots&\bar{\lambda}_{2n}\\
      \hline
      \lambda^{1}&\bar{\lambda}^{1}&&&\\
      \vdots&\vdots&&a^i_j&\\
      \lambda^{2n}&\bar{\lambda}^{2n}&&&
    \end{array}\right)
\end{gather}
where $\bar\lambda_i=J_{\ \ i}^\i\lambda_\i$. The matrix
$L(a,\lambda,\mu)$ is a well-defined $(1,1)$-tensor field on
$\widehat{M}$ (in the sense that after a local
coordinate change in $M$ the components of the matrix $L$ transform
according to tensor rules).

\begin{rem}
  We consider the metric  $g_{ij}$ as a solution of the
  system~\eqref{eq:system} with $\lambda_i=0$ and $\mu=1$. Thus
  \begin{gather}
    L(g,0,1)=\left(
      \begin{array}{cc|ccc}
        1&0&0&\dots&0\\
        0&1&0&\dots&0\\
        \hline
        0&0&&&\\
        \vdots&\vdots&&\delta^i_j&\\
        0&0&&&
      \end{array}\right)={\bf 1}
  \end{gather}
\end{rem}
\begin{rem}
  We see that the matrix $L$ contains as much information as the
  triple $(a_{ij}, \lambda_i, \mu)$, so in a certain sense it is an
  alternative equivalent way to write down the triple. In the next section,
  we will see that the matrix formalism does have advantages: we will
  show that the polynomials of the matrix $L$ also correspond to
  certain solutions of the extended system.

  \par Let us also note  that there is a visually similar construction
  in the theory of projectively equivalent metrics, which uses
  \emph{cone manifolds},
  see~\cite{Matveev2010,Mounoud2010,Cortes2009}. However, in the case
  of $h$-projectively equivalent metrics, the extended operator is not
  covariantly constant (as in the theory of projectively equivalent
  metrics) which poses additional difficulties.
\end{rem}

\subsection{Algebraic properties of $L$}

A linear combination of two matrices of the form \eqref{eq:operator}
is also a matrix of this form, and corresponds to the linear
combination of the solutions (with the same coefficients).  The next lemma shows that the $k$-th
power of the matrix also corresponds to a solution of the extended
system.

\begin{lem}
  \label{lem:product}
  Let $(a, \lambda, \mu)$ be a solution of \eqref{eq:system}. Then,
  for every $k\geq 0$ there exists a solution $(\tilde
  a,\tilde\lambda,\tilde \mu)$ such that
  $$L^k(a,\lambda,\mu)=L(\tilde a,\tilde\lambda,\tilde \mu), \mbox{where }L^k = \underbrace{L
    \cdot ... \cdot L}_\text{$k$ times}.$$
\end{lem}

\begin{proof}
  Given two solutions $(a,\lambda,\mu)$ and $(A,\Lambda,\Mu )$
  of~(\ref{eq:system}), let us calculate the product of the
  corresponding matrices $L(a,\lambda,\mu)$ and $L(A,\Lambda,\Mu )$:
  by direct calculations we obtain
 \begin{multline}
 \label{eq:product}
 L(a,\lambda,\mu)\cdot L(A,\Lambda,\Mu )=\\
 =\left(
 \begin{array}{cc|ccc}
 \mu \Mu +\lambda_k \Lambda^k&\lambda_k \bar\Lambda^k&\mu
 \Lambda_{1}+\lambda_k A^k_1&\dots&\mu\Lambda_{2n}+\lambda_k
 A^k_{2n}\\
 \bar \lambda_k \Lambda^k&\mu \Mu +\lambda_k \Lambda^k&\mu
 \bar \Lambda_{1}+\bar \lambda_k
 A^k_1&\dots&\mu\bar\Lambda_{2n}+\bar\lambda_k A^k_{2n}\\
 \hline
 \vphantom{\dfrac{1}{2}}\Mu \lambda^{1}+a^1_k\Lambda^k&\Mu \bar{\lambda}^{1}+a^1_k\bar\Lambda^k&&&\\
 \vdots&\vdots&\multicolumn{3}{c}{a^i_k
 A^k_j+\lambda^i\Lambda_j+\bar\lambda^i\bar\Lambda_j}\\
 \Mu \lambda^{2n}+a^{2n}_k\Lambda^k&
 \Mu \bar{\lambda}^{2n}+a^{2n}_k\bar\Lambda^k&&&
 \end{array}\right)
 \end{multline}
 Suppose that
 \begin{gather}
 \label{eq:op_eq}
 \mu\Lambda_j+\lambda_k A^k_j= \Mu \lambda_j+a^k_j\Lambda_k\quad\mbox{and}\quad\lambda^k \bar\Lambda_k=0
 \end{gather}
 then
 \begin{gather}
 \label{eq:product_sol}
 L(a,\lambda,\mu)\cdot L(A,\Lambda,\Mu )= L(\underbrace{a^i_k
 A^k_j+\lambda^i\Lambda_j+\bar\lambda^i\bar\Lambda_j}_{\tilde a_{ij}},\underbrace{\mu\Lambda_i+\lambda_k
 A^k_i}_{\tilde \lambda_i},\underbrace{\mu \Mu +\lambda_k \Lambda^k}_{\tilde\mu})
 \end{gather}
 Now we show that the operator $L(\tilde a,\tilde \lambda,\tilde \mu)$ is
 self-adjoint and $\tilde a$, $\tilde \lambda$ and $\tilde \mu$
 satisfy~(\ref{eq:system}).
 \par Indeed, let us check the first equation of~\eqref{eq:system}:
 \begin{multline}
 \tilde a_{ij,k}=(a_{is}
 A^s_j+\lambda_i\Lambda_j+\bar\lambda_i\bar\Lambda_j)_{,k}=
 a_{is,k}A^s_j+a^s_i
 A_{sj,k}+\lambda_{i,k}\Lambda_j+\lambda_i\Lambda_{j,k}+
 \bar\lambda_{i,k}\bar\Lambda_j+\bar\lambda_i\bar\Lambda_{j,k}
 \stackrel{(\ref{eq:system})}{=}\\
 =A^s_j \lambda_i g_{sk}+A^s_j \lambda_s g_{ik}+
 A^s_j \bar \lambda_i J^{s'}_{\ \ s}g_{s' k}+A^s_j \bar\lambda_s
 J^{\i}_{\ \ i}g_{\i k}+\\+a^s_i \Lambda_j g_{sk}+a^s_i \Lambda_s g_{jk}+
 a^s_i \bar \Lambda_j J^{s'}_{\ \ s}g_{s' k}+a^s_i \bar\lambda_s
 J^{\j}_{\ \ j}g_{\j k}+
 \\+
 \mu g_{ik} \Lambda_j-a_{ik}\Lambda_j+\Mu
 g_{jk}\lambda_i-A_{jk}\lambda_i+
 \mu J^{\i}_{\ \ i} g_{\i k} \bar\Lambda_j-J^{\i}_{\ \ i}
 a_{\i k}\bar\Lambda_j+\Mu J^{\j}_{\ \
 j}g_{\j k}\bar\lambda_i-J^{\j}_{\ \ j}A_{\j k}\bar\lambda_i=\\
 =g_{ik}(\lambda_s A^s_j+\mu \Lambda_j)+g_{jk}(\Lambda_s
 a^s_i+\Mu\lambda_i)+
 J^{\i}_{\ \ i}g_{\i k}(\bar\lambda_s A^s_j+\mu \bar\Lambda_j)+
 J^{\j}_{\ \ j}g_{\j k}(\bar\Lambda_s a^s_i+\Mu\bar \lambda_i)
 \stackrel{\eqref{eq:op_eq}}{=}\\
 =\Jij(\tilde\lambda_\i g_{\j k} + \tilde\lambda_\j g_{\i k})
 \end{multline}
 For the second equation one can calculate:
 \begin{multline}
 \tilde \lambda_{i,k}=(\mu\Lambda_i+\lambda_j A^j_i)_{,
 k}=\mu_{,k}\Lambda_i+\mu \Lambda_{i, k}+ \lambda_{j,k}
 A^j_i+\lambda^j A_{ij,k}\stackrel{(\ref{eq:system})}{=}\\
 =-2\lambda_k \Lambda_i+\mu\Mu g_{ik}-\mu A_{ik}+\mu A_{ik}
 -a_{jk}A^j_i+ \lambda^j \Lambda_i g_{j k} + \lambda^j \Lambda_j
 g_{i k}+\lambda^j J^\j_{\ \ j}J^\i_{\ \ i}\Lambda_\i g_{\j k} +
 \lambda^j J^\j_{\ \ j}J^\i_{\ \ i}\Lambda_\j g_{\i k}=\\
 =(\mu \Mu+\lambda^j\Lambda_j)g_{ik}-(\lambda_k \Lambda_i
 +\bar\lambda_k \bar\Lambda_i+A_{ij}
 a^j_k)+\lambda^j\bar\Lambda_j J^\i_{\ \ i}g_{\i
 k}\stackrel{\eqref{eq:op_eq}}{=}
 \tilde \mu g_{ki}-\tilde a_{ki}
 \end{multline}
 From this equation we see that $\tilde a_{ij}$ is symmetric as a
 linear combination of two symmetric tensors. The last equation
 of~\eqref{eq:system} reads
 \begin{multline}
 \tilde \mu_{, i}=(\mu \Mu +\lambda_k \Lambda^k)_{, i}=
 \mu_{, i} \Mu + \mu \Mu_{, i} +\lambda_{k, i} \Lambda^k+\lambda^k
 \Lambda_{k, i}=\\
 \stackrel{(\ref{eq:system})}{=}-2\lambda_i \Mu
 -2\Lambda_i \mu +\Lambda^k (\mu g_{ik}-a_{ik})+\lambda^k (\Mu g_{ik}-A_{ik})=-(\mu\Lambda_i+\lambda_k
 A^k_i)-(\Mu\lambda_i+\Lambda_k
 a^k_i)\stackrel{\eqref{eq:op_eq}}{=}-2\tilde \lambda_i
 \end{multline}
 Thus, $(\tilde a,\tilde \lambda,\tilde \mu)$ is a solution of~\eqref{eq:system}.
 %% Let $L^{k-1}$ satisfy~(\ref{eq:op_eq}), then
 %% $L^{k-1}=L(A,\Lambda,\Mu )$.
 \par Let us now show that the operator
 $L(A,\Lambda,\Mu)=L^k(a,\lambda,\mu)$ satisfies the
 conditions~\eqref{eq:op_eq}.
 \par Since $L^{k}\cdot L=L\cdot L^{k}$, using~(\ref{eq:product}) we obtain
 $$\mu\Lambda_j+\lambda_k A^k_j= \Mu \lambda_j+a^k_j\Lambda_k$$
The last condition will be checked by induction. Suppose $\lambda^i \bar \Lambda_i=0$ then
 $$\lambda^i J^\i_{\ \ i}\widetilde\Lambda_\i=\lambda^i \cdot J^\i_{\
 \ i}(\mu\Lambda_\i+\lambda_k A^k_\i)=\mu\cdot 0 + \lambda^k
 (J_{\ \ i}^\i A_{k \i})\lambda^i=0$$
which completes the proof of Lemma \ref{lem:product}.
\end{proof}
From Lemma \ref{lem:product}, we immediately obtain
\begin{cor}
  Let $(a_{ij },\lambda_i,\mu)$ be a solution of~(\ref{eq:system}) and
  $P(t)=c_k t^k+\dots+c_0$ be an arbitrary polynomial with real
  coefficients. Then there exists a solution $ (A_{ij},\Lambda_i,\Mu)$
  of~(\ref{eq:system}) such that $$L(A_{ij},\Lambda_i,\Mu)= c_k \cdot
  L^k(a_{ij},\lambda_i, \mu) + \dots + \Id:= P(L(a_{ij},\lambda_i,
  \mu)) ,$$ where $\Id$ is the identity $(2n+2)\times
  (2n+2)-$matrix.
\end{cor}

\subsection{There exists a solution $(\check a_{ij}, \check \lambda_{i}, \check \mu)$ such that $L(\check a_{ij}, \check \lambda_i, \check \mu)$ is a projector.}

We assume   that $(M^{2n\geq 4}, g, J)$ is a closed connected Kähler manifold.  Our goal is to show   that 
the existence of a solution $(a_{ij}, \lambda_{i},  \mu)$ of (\ref{eq:system})  such that $\lambda_i\not\equiv 0$ implies the existence of a solution $(\check a_{ij},\check \lambda_i,\check \mu)$ of 
(\ref{eq:system}) such that 
the matrix  $L(\check a_{ij},\check \lambda_i,\check \mu)$ is a non-trivial (i.e. $\neq 0$ and $\neq \Id$) projector. (Recall that a matrix $L$ is a \emph{projector}, if $L^2=L$.).
We need 
\begin{lem}
  \label{lem:minim}
  Let $(M^{2n}, g, J)$ be a connected Kähler manifold and
  $(a_{ij},\lambda_{i},\mu)$  be a solution of \eqref{eq:system}.
  \par Let $P(t)$ be the minimal polynomial of $L(a,\lambda,\mu)$ at
  the point $\hat p \in\widehat M$. Then, $P(t)$ is the minimal
  polynomial of $L(a,\lambda,\mu)$ at every $\hat q\in \widehat M$.
\end{lem}
\begin{con}
We will always assume that the leading coefficient of a minimal polynomial is $1$.
\end{con} 
\begin{proof}
  As we have already proved, there exists a solution $(\tilde a_{ij},\tilde
  \lambda_{i},\tilde \mu)$ such that $$P(L(a,\lambda,\mu))=L(\tilde a,\tilde
  \lambda,\tilde \mu).$$ Since $P(L(a,\lambda,\mu))$ vanishes at the
  point $\hat p=(x_+, x_-,p)$, then  $\tilde a=0$, $\tilde\lambda=0$
  and $\tilde \mu=0$ at   $p$. Then, by
  Lemma~\ref{lem:frobenius}, the solution $(\tilde a_{ij},\tilde\lambda_{i},\tilde
  \mu)$ is identically zero on $M$. Thus, $P(L(a,\lambda,\mu))$
  vanishes at all points of $\widehat M$.  It follows, that the polynomial $P(t)$
  is divisible by the minimal polynomial $Q(t)$  of $L(a,\lambda,\mu)$ at $\hat
  q$. By the same reasoning (interchanging $\hat p$ and $\hat q$), we obtain
  that $Q(t)$ is divisible by $P(t)$.  Consequently, $P(t)=  Q(t)$.
\end{proof}

\begin{cor}
  The eigenvalues of $L(a, \lambda, \mu)$ are constant functions on
  $\widehat M$.
\end{cor}
\begin{proof} By Lemma~\ref{lem:minim}, the minimal  polynomial 
does not depend on  the point of $\widehat M$. Then, the roots of the minimal  polynomial
  are also constant (i.e., do not depend on the point of  $\widehat M$). 
\end{proof}
In order to find the desired special solution of the system (\ref{eq:system}), we will use  that $M$ is closed.
\begin{lem}
  \label{lem:two_eigenvalues}
  Suppose $(M^{2n},g,J)$ is a closed connected Kähler manifold. Let $(a_{ij}, \lambda_i, \mu)$
  be a solution of \eqref{eq:system} such that $\lambda_i\ne 0$ at least at one point.  Then, at every point of $\widehat M$ the matrix $L(a,
  \lambda, \mu)$ has at least two different real eigenvalues.
\end{lem}

\begin{proof}
  Since $M$ is closed, the function $\mu$ admits its maximal  and
  minimal values $\mu_{\max}$ and $\mu_{\min}$.  Let $p \in M$ be a 
  point where $\mu=\mu_{\max}$. At this point, $\mu_{,i}=0$ implying
  $\lambda_i=\bar \lambda_i=0$ in view of the third equation
  of~(\ref{eq:system}). Then, the matrix of $L(a, \lambda, \mu)$ at
  $p$ has the form

  \begin{gather}
    L(a,\lambda,\mu)=\left(
      \begin{array}{cc|ccc}
        \mu_{\max}&0&0&\dots&0\\
        0&\mu_{\max}&0&\dots&0\\
        \hline
        0&0&&&\\
        \vdots&\vdots&&a^i_j&\\
        0&0&&&
      \end{array}\right)
  \end{gather}

  Thus, $\mu_{\max}$ is an eigenvalue of $L(a, \lambda, \mu)$ at $p$
  and, since the eigenvalues are constant, $\mu_{\mbox{\tiny max}}$ is an eigenvalue of $L(a, \lambda, \mu)$ at every point of   $M$. The same
  holds for $\mu_{\min}$. Since $\lambda_i\not\equiv 0$, $\mu$ is not
  constant implying $\mu_{\max}\ne \mu_{\min}$. Finally, $L(a,
  \lambda, \mu)$ has two different real eigenvalues $\mu_{\max},
  \mu_{\min}$ at every point.
\end{proof}

\begin{rem}
  \label{-3}
  For further use let us note that 
   in the proof of Lemma~\ref{lem:two_eigenvalues} we have
  proved that if $\mu_{,i}=0$ at a point $p$ then $\mu(p)$ is an
  eigenvalue of $L$.
\end{rem}
Finally, let us show that there is always a solution of (\ref{eq:system}) of the desired special kind:
\begin{lem}
  Suppose $(M^{2n},g,J)$ is a closed and connected Kähler manifold. For every solution $(a_{ij}, \lambda_i, \mu)$ of \eqref{eq:system}
  such that $\lambda_i$ is not identically zero on $M$, there exists a polynomial $P(t)$
  such that $P(L(a, \lambda, \mu))$ is a non-trivial (i.e. it is neither
  ${\bf 0}$ nor $\Id$) projector.
\end{lem}

\begin{proof}
  We take a point $\hat{p}\in \widehat{M}$. By Lemma~\ref{lem:two_eigenvalues},
  $L(a_{ij}, \lambda_i, \mu)$ has at least two real eigenvalues at the
  point $\hat{p}$.  Then, by linear algebra, there exists a polynomial $P$
  such that $P( L(a_{ij}, \lambda_i, \mu))$ is a nontrivial projector
  at the point $p$. Evidently, a matrix $C$ is a nontrivial projector,
  if and only if its minimal polynomial is $t(t-1)$ (multiplied by any
  nonzero constant). Since by Lemma~\ref{lem:minim} the minimal polynomial
  of $P( L(a_{ij}, \lambda_i, \mu))$ is the same at all points, the
  matrix $P( L(a_{ij}, \lambda_i, \mu))$ is a projector at every point
  of $M$.
\end{proof}

Thus, (under the assumptions of Theorem~\ref{pos}), without loss of
generality we can think that a solution of the
system~(\ref{eq:system}) on a closed and connected Kähler manifold $M$ with degree
of mobility $\geq 3$ is chosen  such that the corresponding $L$ is a
projector.

\subsection{Structure of eigenspaces of $a_j^i$, if $L(a, \lambda, \mu)$ is a  nontrivial projector}

We assume that    $L(a, \lambda, \mu)$ is  a nontrivial  projector. 
Then, it has precisely  two
eigenvalues: $1$ and $0$ and the $(2n+2)$-dimensional tangent space of
$\widehat M$ at every point $\hat x=(x_+,x_-,p)$ can be decomposed
into the sum of the corresponding eigenspaces $$T_{\hat x}\widehat M =
E_{L(a, \lambda, \mu)}(1)\oplus E_{L(a, \lambda, \mu)}(0).$$ 

The dimensions of $E_{L(a, \lambda, \mu)}(1)$ and of $E_{L(a, \lambda, \mu)}(0)$ are even; we assume that the dimension  of $E_{L(a, \lambda, \mu)}(1)$ is $2k+2$ and the dimension of $E_{L(a, \lambda, \mu)}(0)$ is $2n-2k$. 

By Lemma~\ref{lem:two_eigenvalues}, $\mu_{\max}$ and $\mu_{\min}$ are
eigenvalues of $L(a, \lambda, \mu)$. Then, $\mu_{\min}=0\leq\mu(x)\leq
1=\mu_{\max}$ on $M$. In view of Remark~\ref{-3}, the only critical
values of $\mu$ are $1$ and $0$.

\begin{lem}
Let $(a_{ij},\lambda_{i},\mu)$ be a solution of (\ref{eq:system}) such that $L(a,\lambda,\mu)$ is a non-trivial projector.   
Then, the following statements hold:
  \label{lem:eigen}
  \begin{enumerate}
  \item At a point $p$ such that $0<\mu< 1$,
    $a^i_j$ has the following structure of eigenvalues and eigenspaces
    \begin{enumerate}
    \item eigenvalue $1$ with geometric multiplicity $2k$;% $E_1=E_{L(a, \lambda, \mu)}(1)\cap TM,$
    \item eigenvalue $0$ with geometric  multiplicity $(2n-2k-2)$;% $E_0=E_{L(a, \lambda, \mu)}(0)\cap TM$,
    \item eigenvalue $(1-\mu)$ with multiplicity $2$.% $E_{1-\mu}={\mathrm span}\{\lambda^i,\bar\lambda^i\}$.
    \end{enumerate}
  \item At a  point $p$ such that $\mu=1$, $a^i_j$ has the following
    structure  of eigenvalues and eigenspaces:
    \begin{enumerate}
    \item eigenvalue $1$ with geometric multiplicity $2k$;% $E_1=E_{L(a, \lambda, \mu)}(1)\cap TM,$
    \item eigenvalue $0$ with geometric  multiplicity $(2n-2k)$;% $E_0=E_{L(a, \lambda, \mu)}(0)\cap TM$.
    \end{enumerate}
  \item At a  point $p$ such that $\mu=0$, $a^i_j$ has the following
    structure of   eigenvalues and eigenspaces:
    \begin{enumerate}
    \item eigenvalue $1$ with geometric multiplicity $2k+2$;
    \item eigenvalue $0$ with geometric multiplicity $(2n-2k-2)$.
    \end{enumerate}
  \end{enumerate}
\weg{  (see Fig.~\ref{fig:eigenspaces}).}
\end{lem}

\begin{con} We identify  $M$  with the set $(0,0)\times M\subset \widehat M$. 
This identification allows us to consider   $T_xM$ as  a  linear subspace of $T_{(0,0)\times x} \widehat M$: the  vector $(v_1,...,v_n)\in T_xM$ is identified with $(0,0,v_1,...,v_n)\in T_{(0,0)\times x} \widehat M$. 
\end{con} 

\begin{proof}
  For any vector $v\in E_1=E_{L(a, \lambda, \mu)}(1)\cap TM$ we  calculate
  \begin{gather}
    L(a, \lambda, \mu) v=\left(
      \begin{array}{cc|ccc}
        \mu&0&\lambda_{1}&\dots&\lambda_{2n}\\
        0&\mu&\bar{\lambda}_{1}&\dots&\bar{\lambda}_{2n}\\
        \hline
        \lambda^{1}&\bar{\lambda}^{1}&&&\\
        \vdots&\vdots&&a^i_j&\\
        \lambda^{2n}&\bar{\lambda}^{2n}&&&
      \end{array}\right)
    \left(
      \begin{array}{c}
        0\\
        0\\
        v^1\\
        \vdots\\
        v^{2n}
      \end{array}
    \right)=
    \left(
      \begin{array}{c}
        \lambda_j v^j\\
        \bar \lambda_j v^j\\
        \\
        a^i_j v^j\\
        \\
      \end{array}
    \right)
    = \left(
      \begin{array}{c}
        0\\
        0\\
        v^1\\
        \vdots\\
        v^{2n}
      \end{array}
    \right)
  \end{gather}

  Thus, $v=(v^1,\dots v^{2n})$ is an eigenvector of $a^i_j$ with
  eigenvalue $1$. Moreover, it is orthogonal to both $\lambda^i$ and
  $\bar\lambda^i$.   Similarly, any $v \in E_0=E_{L(a, \lambda,
    \mu)}(0)\cap T_xM$ is an eigenvector of $a^i_j$ with eigenvalue~$0$
  and it is orthogonal to $\lambda^i$ and $\bar\lambda^i$. Note that the dimension of $E_1$ is at least  
  $\dim E_{L(a, \lambda, \mu)}(1) -  2 = 2k$, and the dimension of $E_0$ is at least $\dim E_{L(a, \lambda, \mu)}(0) -  2 = 2n - 2k - 2$.

  Thus,  at every point $x$  there are three pairwise orthogonal subspaces in $T_xM$: $E_1$,
  $E_0$ and $span\{\lambda^i,\bar\lambda^i\}$.

If  $0<\mu<1$ at $x$, $\lambda_i\ne 0$  by Remark \ref{-3}. Then, the dimension of $E_1\oplus E_0\oplus \spann\{\lambda^i,\bar \lambda^i\}$   is at least  $2n - 2k - 2 + 2k +2 = 2n$.  Since $E_1\oplus E_0\oplus \spann\{\lambda^i,\bar \lambda^i\}\subseteq  T_{ x} M,$  
the dimension of $E_1$ is $  2n-2k-2$  and the dimension of $E_0$ is $2k$, and $E_1\oplus E_0\oplus \spann\{\lambda^i,\bar \lambda^i\}=  T_{x} M$.

 Let us now show that 
   $\lambda^i$ and $\bar\lambda^i$ are eigenvectors of
  $a^i_j$ with the eigenvalue $(1-\mu)$. We 
   multiply  the first basis vector $(1,0,\dots,0)$ by 
    the matrix $L(a, \lambda,
  \mu)^2-L(a, \lambda, \mu)$ (which is identically zero). We obtain 
  \begin{gather}
    0 = (L(a, \lambda, \mu)^2-L(a, \lambda, \mu)) \left(
      \begin{array}{c}
        1\\
        0\\
        0\\
        \vdots\\
        0
      \end{array}
    \right)=
    \left(
      \begin{array}{c}
        \mu^2 + \lambda_j \lambda^j - \mu\\
        \bar\lambda_i\lambda^i\\
        \\
        \mu \lambda^i + a^i_j \lambda^j-\lambda^i\\
        \\
      \end{array}
    \right)
  \end{gather}
  This gives us the necessary equation $a^i_j \lambda^j= (1-\mu)\lambda^i$.

Finally, we have that $T_xM$ is the direct sum $E_1\oplus E_0\oplus \spann\{\lambda^i,\bar \lambda^i\}; $
 $E_1$  consists of eigenvectors of $a^i_j$ with eigenvalue $1$ and has dimension 
 $2n - 2k - 2$,  $E_1$  consists of eigenvectors of $a^i_j$ with eigenvalue $0$  and has dimension  $2k$; 
 $\spann\{\lambda^i,\bar \lambda^i\}$ consists of eigenvectors of $a^i_j$ with eigenvalue $(1-\mu)$ and has dimension $2$, as we claimed in the first statement of the lemma. 
 
 The proof at the points $x$ such that  $\mu(x)=0$ or  $\mu(x)=1$ is similar (and is easier), and will be left to the reader.
\end{proof}

\subsection{ If there exists a solution $(a,\lambda,\mu)$ of
 the system~(\ref{eq:system}) corresponding to a non-trivial projector,
 the metric $g$ is positively definite on $M$ (assumed closed).}

Above we have proved that, under the assumptions of Theorem \ref{pos}, there always exists a solution $(a_{ij},\lambda_{i},\mu)$ of (\ref{eq:system})
such that the corresponding matrix $L(a,\lambda,\mu)$ is a non-trivial
projector, implying that the eigenvalues and the dimension of eigenspaces of $a^i_j$ is given by Lemma \ref{lem:eigen}. Now we are ready to prove that $g$ is positively definite
(as we claimed in Theorem~\ref{pos}).

Let us consider such a solution $(a_{ij},\lambda_{i},\mu)$. We rewrite the second equation in (\ref{eq:system}) in the form
\begin{gather}
  \label{eq:mu}
  \mu_{,i j} =2 a_{i j} -2 \mu\, g_{i j}
\end{gather}
Let $p$ be a point where $\mu$ takes its maximal  value $1$. As we have already shown, 
$\lambda^i(p)=0$ and the tangent space $T_pM$ is  the direct sum of
the eigenspaces of $a^{i}_j$:

$$T_p M=E_1\oplus E_0$$
Consider the restriction of~(\ref{eq:mu}) to $E_0$. Since the restriction of
the bilinear form $a_{i j}$ to $E_0$ is identically zero, the restriction of~(\ref{eq:mu}) to $E_0$
reads
$$\left.\mu_{,i j}\right|_{E_0}=-2 \left.g_{i j}\right|_{E_0}.$$

Now,  $\mu_{,i j}$ is the Hessian of $\mu$ at the maximum point
$p$. Then, it is non-positively definite. Hence, the non-degenerate metric
tensor $g_{ij}$ is positively definite on $E_0$ at $p$.
Let us now consider the distribution of the orthogonal complement
$E_1^\bot$, which is well-defined, smooth and integrable on $\{x\in M\mid \mu(x)>0\}$. 
 The restriction of the metric $g$ to  $E_1^\bot$ is
non-degenerate at  the points of  $\{x\in M\mid \mu(x)>0\}$.  Since at the  point $p$ 
  $E_1^\bot$ coincides with $E_0$, it  is positively definite at $p$.  Hence, by
continuity,  it is positively definite at   the connected component  of 
 $\{x\in M\mid \mu(x)>0\}$  containing $p$. Since every connected component of  $\{x\in M\mid \mu(x)>0\}$ has a point such that $\mu=1$, the restriction of the metric $g$ to  $E_1^\bot$ is
positively definite  at all points of   $\{x\in M\mid \mu(x)>0\}$. 

Similarly, at a minimum point $q$ one can consider the restriction
of~(\ref{eq:mu}) to $E_1$:
$$\left.\mu_{,i j}\right|_{E_1}=2 \left.a_{i j}\right|_{E_1}=2
\left.g_{i j}\right|_{E_1},$$
since $a^{i}_{j\,|E_{1}}=\delta^{i}_{j\,|E_{1}}$.  Then,  $g$ is positively definite on $E_1$ at $q$.
Considering the distribution $E_0^\bot$, we obtain that the restriction of $g$
to $E_0^\bot$ is positively definite at   
 $\{x\in M\mid \mu(x)<1\}$. 
 
Evidently, the   sets $\{x\in M\mid \mu(x)<1\}$  and $\{x\in M\mid \mu(x)>0\}$  have  an nonempty intersection. 
At every point $x$ of the intersection, $T_xM= E_0^\bot + E_1^\bot$. Since the restriction of the metric to $E_0^\bot$ and to  $E_1^\bot$ is positively definite, the metric is positively definite as we claimed. 
 Theorem~\ref{pos} is proved.

\section{Tanno-Theorem completes  the proof of Theorem 1}
\label{sec:tan}

We assume that $(M^{2n}, g, J)$ is a closed connected Kähler manifold of
dimension $2n \geq 4$ with degree of mobility $D \geq 3$.  Let $\bar g$
be a metric $h$-projectively equivalent to $g$. We consider the
corresponding solution $(a_{ij}, \lambda_i, \mu)$ of the extended
system. If  the metric $\bar g$ is  not  affinely equivalent
to $g$, by Theorem~\ref{thm:B=0} we obtain $B\ne 0$. As we
explained in the beginning of Section~\ref{sec:proj}, by multiplication of the metric by a nonzero constant, we can achieve  $B=-1$. Without loss of generality, we think that $B=-1$. 
 By Theorem~\ref{pos}, the metric $g$ is Riemannian.

Now, by Lemma~\ref{lem:tanno}, the function $\lambda:=\tfrac{1}{4} a^i_i$
satisfies the equation
\begin{align}
  \lambda_{,ijk}+(2\lambda_{,k}g_{ij}+\lambda_{,i}g_{jk}+\lambda_{,j}g_{ik}
  +(J_{\ \ i}^{\alpha}J_{\ \ j}^{\beta}+J_{\ \ j}^{\alpha}J_{\ \ i}^{\beta})\lambda_{,\alpha}g_{\beta k})=0,
\end{align}
moreover, by Remark~\ref{aff},  if $\bg$ is not affinely equivalent to
$g$, the function $\lambda$ is not a constant.

As we recalled in Section \ref{1.8}, this equation was considered in~\cite{Tanno1978}. Tanno has
proved, that the existence of a non-constant solution of this equation
on a closed connected Riemannian manifold implies that the metric $g$
has positive constant holomorphic sectional curvature  equal to 4 (see~\cite[Theorem
10.5]{Tanno1978}, and also Section \ref{1.8}). Then, $(M^{2n}, g, J)$ is 
$(\mathbb{C}P(n), 4 \cdot g_{FS} , J_{standard})$ as we
claimed. Theorem \ref{thm:main} is proved.

\begin{rem} \label{tmp}  As we already mentioned in Section \ref{plan}, in Sections \ref{sec:B=0}, \ref{sec:proj} we did not actually use the assumption that the degree of mobility is $\ge 3$: we  used the system \eqref{eq:mg+Ba} only. Thus, the following statement holds: 

{\it  Let $(M^{2n\ge 4}, g, J)$ be a closed connected Kähler manifold. Assume there exists a solution $(a_{ij}, \lambda_i, \mu)$ of \eqref{eq:mg+Ba} such that $\lambda_i\not\equiv 0$. Then, $(M^{2n}, g, J)$ is  $(\mathbb{C}P(n), \const\cdot  g_{FS}, J_{standard})$ (for a certain $\const \ne 0$).

}

\end{rem} 

\section{Proof of Theorem~\ref{thm:tanno}: equation~\eqref{eq:tanno} is equivalent to  system~\eqref{eq:mg+Ba}} \label{newtan} 

In Lemma~\ref{lem:tanno}, we have shown that for a solution of the
extended system~\eqref{eq:mg+Ba}  equation~\eqref{eq:tanno} is
fulfilled. We will now show that  a nonconstant solution of~\eqref{eq:tanno} allows us 
to construct a solution $(a_{ij},\lambda_i,\mu)$ of the extended
system~\eqref{eq:mg+Ba} (with $B=\kappa$ and $\lambda_i=f_{,i}\not\equiv 0$) 
 provided that the manifold is closed.

 Let $f$ be a non-constant solution of equation~\eqref{eq:tanno}
on a closed connected  manifold $M$. Then,  $\kappa\ne 0$. Indeed, we can proceed as in Section \ref{sec:B=0}: if $\kappa=0$, then equation \eqref{eq:tanno} reads $f_{,ijk}=0$.  Then, the hessian $f_{,ij} $ of the
  function $f$ is covariantly constant.
 Since the manifold is closed, the function $f$ has a
  minimum and a maximum. At a minimum, the Hessian must be
  non-negatively definite, and at a maximum it must be nonpositively
  definite. Therefore the Hessian is null, and $f_{,i}$ is
  covariantly constant. But as it vanishes at the extremal points, it
  vanishes everywhere. Thus, $f = \const $ contradicting the assumptions.

  Consider the 
symmetric, hermitian tensor $a_{ij}$ defined by the following formula:
\begin{gather} \label{-01}
  a_{ij}=\frac{1}{\kappa}f_{,ij}-2  f g_{ij}
\end{gather}
Let us check that $(a_{ij}, \lambda_i= f_{,i}, \mu= 2  \kappa f)$ satisfies~\eqref{eq:mg+Ba} with $B=\kappa$.  Indeed,  covariantly differentiating $a_{ij}=\frac{1}{\kappa}f_{,ij}-2 f g_{ij}$ and substituting 
\eqref{eq:tanno}, we obtain  
\begin{multline}
  a_{ij,k}=\frac{1}{\kappa}f_{,ijk}-2 f_{,k} g_{ij}=2 f_{,k} g_{ij}+ f_{,i} g_{j k} +
    f_{,j} g_{i k} - \bar{f}_{,i}J_{jk}-\bar{f}_{,j}J_{ik} -2 f_{,k}
    g_{ij}=\\
    =f_{,i} g_{j k}+f_{,j} g_{i k}-\bar{f}_{,i}J_{jk}-\bar{f}_{,j}J_{ik},
\end{multline} which is the first
equation of  \eqref{eq:mg+Ba}. 
The second   equation  of  \eqref{eq:mg+Ba}  is   equivalent to  \eqref{-01}, the third equation is 
fulfilled by the construction. 
 Since $f$ is non-constant,
$\lambda_i=f_{,i}\not \equiv  0$.  Now, as we proved in  Section \ref{sec:proj}, the metric $-\textrm{sgn}(B) \cdot g$ is positively definite. Finally, for positively definite metrics,  Theorem~\ref{thm:tanno} was proved by Tanno  in \cite{Tanno1978}.      Theorem~\ref{thm:tanno} is proved.

\subsection*{\bf Acknowledgements.} We thank D.V. Alekseevsky, B. Kruglikov, D. Calderbank, V. Cortes, K. Kiyohara and P. Topalov for usefull discussions and the anonymous referee for his valuable suggestions.

\nocite{*}
\bibliographystyle{plain}

\end{document}